\documentclass[fleqn,11pt,twoside]{article}

\usepackage{amsthm,amsthm,amssymb, color, xcolor,epsfig, graphics, subfigure}
\usepackage{amsmath, graphicx, latexsym, lscape}

\makeatletter
\newcommand{\copyrightnote}[2]{{\renewcommand{\thefootnote}{}
 \footnotetext{\small\it
\begin{flushleft}
 \copyright \ #1   #2
\end{flushleft}}}}

\newcommand{\Name}[1]{\begin{flushleft}
                       \LARGE \bf #1
                       \end{flushleft}\vspace{-3mm}}

\newcommand{\Author}[1]{\begin{flushleft}
                       \it #1 \end{flushleft}}

\newcommand{\Address}[1]{\begin{flushleft}
                       \it #1 \end{flushleft}}

\newcommand{\Date}[1]{\begin{flushleft}
                      \small  \it #1 \end{flushleft}}

\newcommand{\evenhead}{Author \ name}
\newcommand{\oddhead}{Article \ name}

\renewcommand{\@evenhead}{
\hspace*{-3pt}\raisebox{-15pt}[\headheight][0pt]{\vbox{\hbox to \textwidth
{\thepage \hfil \evenhead}\vskip4pt \hrule}}}
\renewcommand{\@oddhead}{
\hspace*{-3pt}\raisebox{-15pt}[\headheight][0pt]{\vbox{\hbox to \textwidth
{\oddhead \hfil \thepage}\vskip4pt\hrule}}}
\renewcommand{\@evenfoot}{}
\renewcommand{\@oddfoot}{}

\long\def\@makecaption#1#2{%
  \vskip\abovecaptionskip
  \sbox\@tempboxa{\small \textbf{#1.}\ \ #2}%
  \ifdim \wd\@tempboxa >\hsize
    {\small \textbf{#1.}\ \ #2}\par
  \else
    \global \@minipagefalse
    \hb@xt@\hsize{\hfil\box\@tempboxa\hfil}%
  \fi
  \vskip\belowcaptionskip}

\newcommand{\JNMPnumberwithin}[3][\arabic]{%
  \@ifundefined{c@#2}{\@nocounterr{#2}}{%
    \@ifundefined{c@#3}{\@nocnterr{#3}}{%
      \@addtoreset{#2}{#3}%
      \@xp\xdef\csname the#2\endcsname{%
        \@xp\@nx\csname the#3\endcsname .\@nx#1{#2}}}}%
}

\newcommand{\resetfootnoterule} {
  \renewcommand\footnoterule{%
  \kern-3\p@
  \hrule\@width.4\columnwidth
  \kern2.6\p@}
}

\renewcommand{\footnoterule}{}

\makeatother

\theoremstyle{definition}

\definecolor{Red}{rgb}{0.9,0.2,0.1}
\definecolor{dkRed}{rgb}{0.5,0.2,0.4}
\definecolor{Green}{rgb}{0.2,0.9,0.2}
\definecolor{dkGreen}{rgb}{0.1,0.8,0.1}
\definecolor{Yellow}{rgb}{1,1,0}
\definecolor{Navy}{rgb}{0.1,0.1,0.4}
\definecolor{Navy2}{rgb}{0.1,0.1,0.5}
\definecolor{Black}{rgb}{0,0,0}
\definecolor{Orange}{rgb}{1,0.55,0.02}
\definecolor{Pink}{rgb}{0.88,0.09,0.77}
\definecolor{Grey}{rgb}{0.7,0.7,0.7}
\definecolor{Math}{rgb}{0.07,0.63,0.08}
\definecolor{Violet}{rgb}{0.4,0,0.7}

\begin{document}

\renewcommand{\evenhead}{ {\LARGE\textcolor{blue!10!black!40!green}{{\sf \ \ \ ]ocnmp[}}}\strut\hfill
W Hereman and R Naz
}

\renewcommand{\oddhead}{ {\LARGE\textcolor{blue!10!black!40!green}{{\sf ]ocnmp[}}}
\ \ \ \ \
Conservation laws of nonlinear PDEs from elasticity and acoustics}
%

\thispagestyle{empty}
\newcommand{\FistPageHead}[3]{
\begin{flushleft}
\raisebox{8mm}[0pt][0pt]
{\footnotesize \sf
\parbox{150mm}{{\textcolor{blue!10!black!40!green}{
{\bf Open Communications in Nonlinear Mathematical Physics}}}
\ \ {Special Issue: Bluman}, 2025\\[0.1cm]
\strut\hfill ocnmp: 17124, pp #2\hfill {\sc #3}}}\vspace{-13mm}
\end{flushleft}}

\FistPageHead{1}{\pageref{firstpage}--\pageref{lastpage}}{ \ \ }

\strut\hfill

\strut\hfill

\copyrightnote{The authors. Distributed under a Creative Commons Attribution 4.0 International License}

\begin{center}

{\bf {\large A Special OCNMP Issue in Honour of George W Bluman}}
\end{center}

\smallskip

\Name{Conservation laws of nonlinear PDEs arising in elasticity and acoustics
in Cartesian, cylindrical, and spherical geometries}

\Author{Willy Hereman$^{a1}$, Rehana Naz$^{b}$}

\Address{
$^{a}$Department of Applied Mathematics and Statistics,
Colorado School of Mines, Golden CO 80401-1887, USA
\\
Email: whereman@mines.edu
\\
$^{1}$Corresponding author
\\
$^{b}$Independent Researcher, Pasadena, CA 91107, USA
\\

Email: rehananaz$\_$qau@yahoo.com
}

\Date{Received December 16, 2025; Accepted December 21, 2025}

\setcounter{equation}{0}

\smallskip

\noindent {\bf Citation format for this Article:}
\newline W Hereman and R Naz, Conservation laws of nonlinear PDEs arising in elasticity
and acoustics in Cartesian, cylindrical, and spherical geometries,
{\it Open Commun. Nonlinear Math. Phys.}, Special Issue:\, Bluman,
ocnmp: 17124, \pageref{firstpage}--\pageref{lastpage}, 2025.

\strut\hfill

\noindent {\bf The permanent Digital Object Identifier (DOI) for
this Article:}
\newline 
{\it 10.46298/ocnmp.17124} \strut\hfill

\begin{abstract}
Conservation laws are computed for various nonlinear partial
differential equations that arise in elasticity and acoustics. Using
a scaling-homogeneity approach, conservation laws are established
for two models describing shear wave propagation in a circular
cylinder and a cylindrical annulus. Next, using the multiplier
method, conservation laws are derived for a parameterized system of
constitutive equations in cylindrical coordinates involving a
general expression for the Cauchy stress. Conservation laws for the
Khokhlov-Zabolotskaya-Kuznetsov equation and Westervelt-type
equations in various coordinate systems are also presented.
\end{abstract}

\label{firstpage}

\section{Introduction}
\label{intro}
 One of the beautiful applications of symmetry methods
to partial differential equations (PDEs) is the computation of
conservation laws which has been a central theme in the work of
Prof.\ Bluman and his collaborators
\cite{Bluman-Anco-book-2002,Bluman-Cheviakov-Anco-book-2010,Cheviakov-Zhao-book-2024}.
Conservation laws consist of conserved densities and associated
fluxes. The integral of a conserved density often has a physical
meaning. Depending on the circumstances, it might express
conservation of energy, momentum, angular momentum, and the like.
Conservation laws have many uses in the study of PDEs such as
establishing the existence and uniqueness of solutions and
investigating their stability. They also play a key role in the
development of accurate numerical methods to solve nonlinear PDEs.

Of the many techniques
\cite{Anco-FIC-chapter-2017,Bluman-Cheviakov-Anco-book-2010,Cheviakov-JEM-2010,
Naz-Mahomed-Mason-AMC-2008,Wolf-EJAM-2002}
to compute conservation laws, we focus on two approaches:
(i) a scaling-homogeneity method
\cite{goktas-hereman-jsc-1997,hereman-etal-book-nova-2009, hereman-etal-book-birkhauser-2005,poole-hereman-jsc-2011}
which uses linear algebra and does not require solving PDEs, and
(ii) the multiplier method
\cite{Anco-Bluman-EurApplMath-II-2002,Naz-Mahomed-Mason-AMC-2008,Olver-1993,Steudel-1962}
which is analogous to the integrating factor method \cite{Bluman-Anco-book-2002}
for ordinary differential equations.

Applying these methods which are implemented in the codes
\verb|ConservationLawsMD.m| \cite{hereman-poole-conservationlawsMD}
and \verb|GeM| \cite{Cheviakov-website-2025}, respectively, 
we compute conservation laws of various nonlinear PDEs that arise 
in elasticity and acoustics. 
In \cite{naz-hereman-pdeam-2024} we derived
conservation laws based on a constitutive equation modeling stress
in elastic materials with a geometry suitable for Cartesian
coordinates. Here we compute conservation laws for PDEs that model
the propagation of elastic waves in a circular cylinder and
cylindrical annulus. More precisely, we focus on models studied by
Kambapalli et al.\ \cite{Kambapalli-Kannan-Rajagopal-2014} and Magan
et al.\ \cite{Magan-Mason-Harley-MMS-2020}. These models belong to a
general class of implicit constitutive equations where the strain
$\epsilon$ is expressed as a non-invertible function $F(\sigma)$ of
the stress $\sigma$. These so-called {\emph{implicit theories}}
originated in work of Rajagopal
\cite{rajagopal-AP-2003,rajagopal-ZAMP-2007} and have since been
advocated by many others (see, e.g.,
\cite{bustamante-in-book-2026,bustamante-rajagopal-in-book-2020,
bustamante-sfyris-MMS-2015} and the references therein).

We start with computing conservation laws for the two models
reported in \cite{Kambapalli-Kannan-Rajagopal-2014} and
\cite{Magan-Mason-Harley-MMS-2020}. Both models are formulated in
cylindrical coordinates and assume a power-law constitutive
relation. Next, to capture other constitutive relations reported in
the literature (see Table 2), we replace the power-law by an
arbitrary function, $\epsilon = F(\sigma)$, relating the scalar
linearized strain to Cauchy stress. 
The conserved densities we obtain are independent of $F$ but 
the fluxes depend on the integral of $F$. 
For completeness, we also compute conservation laws
of the single wave equation for the stress that arises upon
elimination of the displacement from the equations of motion.

A couple of equations from nonlinear acoustics are also considered.
The first one is the Khokhlov-Zabolotskaya-Kuznetsov (KZK) equation
\cite{Hamilton-ch8-springer-2024,
Hamilton-Morfey-ch3-springer-2024,Rozanova-Pierrat-CMS-2009}
(sometimes referred to as the two-dimensional Burgers equation)
which models the propagation of sound beams in various nonlinear
media, in particular, sound waves generated by parametric acoustic
areas. Other applications include the propagation of ultrasound in
dissipative media, long waves in ferromagnetic media, and pulsed
sound beams in thermo-viscous media. Conservation laws are computed
for the KZK equation in Cartesian coordinates as well as cylindrical
and spherical coordinates. It turns out that for the KZK equation in
Cartesian coordinates, the conservation laws have a couple of
space-dependent coefficients; one must be a harmonic function while
the other must be a solution of a Poisson equation. Since the
Laplace equation has infinitely many solutions there are infinitely
many conservation laws. In cylindrical and spherical coordinates the
equations for these coefficients are slightly more complicated but
the conclusion remains the same: The KZK equation in plane,
cylindrical, and spherical geometries has infinite conservation
laws.

The second equation is the Westervelt equation
\cite{Hamilton-Morfey-ch3-springer-2024,Kaltenbacher-EECT-2025}
named after Peter Westervelt who received the Nobel Prize in Physics
in 1964 for his contributions to quantum mechanics. His equation is
used to describe the nonlinear propagation of pressure waves in
nonlinear media. It has many engineering applications in underwater
acoustics, sonar antennas, infrasound in the atmosphere, among
others. The equation also models the propagation of high-intensity
ultrasound in tissue with numerous medical applications.

The dissipative version of Westervelt's equation in one spatial
dimension has been the subject of recent studies by Anco et al.\
\cite{Anco-Marquez-Garrido-Gandarias-2023} and M\'{a}rquez et al.\
\cite{Marquez-Recio-Gandarias-2025}. Anco and co-workers not only
report new local and non-local conservation laws but also provide an
in-depth analysis of various types of symmetries including hidden
symmetries for the potential equation, hidden variational
structures, recursion and Noether operators, etc. A similar study is
carried out in \cite{Marquez-Recio-Gandarias-2025} for a slightly
different version of the one-dimensional Westervelt equation with an
arbitrary nonlinearity.

In this paper we compute conservation laws of the Westervelt
equation in more than one space variable. Our computational results
confirm recent findings by Sergyeyev \cite{Sergyeyev-QTDS-2024} who
gives a complete description of local conservation laws of the
multi-dimensional dissipative Westervelt equation by using the
direct method \cite{Anco-Bluman-EurApplMath-I-2002,Olver-1993} and a
theorem from Igonin \cite{Igonin-JPA-2024}. In particular, Sergyeyev
has shown that an infinite number of conservation laws exists for
the case of two or more space variables. This is in contrast with
the one-dimensional dissipative Westervelt equation which has only a
finite number of local conservation laws
\cite{Anco-Marquez-Garrido-Gandarias-2023}.

Using the adjoint-symmetry approach, Anco \cite{Anco-EurApplMath-II-2023}
computed conservation laws of Westervelt's equation in spherical coordinates
with the standard quadratic nonlinearity.
In this paper we get similar results for the Westervelt's equation
(in both spherical and cylindrical coordinates) where the quadratic nonlinear
term has been replaced with an arbitrary nonlinear function.

This paper is organized as follows. In
Section~\ref{model-eqs-cons-laws}, conservation laws are given for
the two forementioned models for shear wave propagation in a
circular cylinder and cylindrical annulus. The conservation laws for
these two models are computed in Section~\ref{cons-laws-scaling}
using the scaling-homogeneity approach. 
In Section~\ref{cons-laws-multiplier}, conservation laws are computed
with the multiplier method for a parameterized system involving an
arbitrary function of stress 
covering a broad class of models for wave propagation in elastic
materials formulated in terms of cylindrical coordinates.
For specific values of the parameters, the general system includes the 
two previously mentioned models.
Conservation laws for the KZK equation in Cartesian, cylindrically,
and spherical coordinates are derived in
Section~\ref{application-KZK}. In
Section~\ref{application-Westervelt}, conservation laws are computed
with the multiplier method for Westervelt-type equations in
multi-space (Cartesian) coordinates as well as cylindrical and
spherical coordinates. Finally, a brief discussion of the results
and some conclusions are given in Section~\ref{conclusions}.

\section{Model equations and some of their conservation laws}
\label{model-eqs-cons-laws} In this section we consider two models
for shear wave propagation in a circular cylinder and cylindrical
annulus. The constitutive equations for these geometries are
appropriately expressed in cylindrical coordinates. The models under
consideration are based on {\emph{implicit}} constitutive equations
where the strain is expressed as a non-invertible function of the
stress. They were introduced by Rajagopal
\cite{rajagopal-AP-2003,rajagopal-ZAMP-2007} and have numerous
applications (see, e.g.,
\cite{bustamante-in-book-2026,bustamante-rajagopal-in-book-2020}).

\subsection{Model due to Kambapalli et al.\
\cite{Kambapalli-Kannan-Rajagopal-2014} with conservation laws}
\label{model-kambapalli-cons-laws}
The first model \cite{Kambapalli-Kannan-Rajagopal-2014} in
non-dimensional variables,
\begin{eqnarray}
\label{kam-(2.18)}
\sigma_r+\frac{2\sigma}{r} &=& \delta u_{tt},
\\
\label{kam-(2.17)}
u_r-\frac{u}{r} &=& \tfrac{1}{\delta} \sigma (\beta + \sigma^2)^n,
\end{eqnarray}
describes shear waves in a cylindrical annular region modeled in
terms of cylindrical coordinates $(r,\theta, z)$ and time $t.$ As
usual, subscripts denote partial derivatives, e.g., $\sigma_r =
\tfrac{\partial \sigma}{\partial r}$ and $u_{tt} = \tfrac{\partial^2
u}{\partial t^2}$.

The first equation expresses the balance of linear momentum. The
second equation is a consequence of the constitutive relation. These
equations of motion for stress $\sigma(r,t)$ and displacement
$u(r,t)$ are for a specific (but rather simple) model where the
linear strain is a nonlinear function of stress expressed as a power
law where the arbitrary exponent $n \ge 0$ is a natural or rational
number (see also
\cite{kannan-rajagopal-saccomandi-wm-2014,magan-etal-wm-2018}).

The reciprocal of constant parameter $\delta$, i.e.,
$\tfrac{1}{\delta} = \frac{\alpha}{\sqrt{\gamma}}$, is the displacement
gradient which involves two material parameters $\alpha \ge 0$ and $\gamma \ge 0$.
An auxiliary constant parameter $\beta$ has been
introduced\footnote{Once the computations are done one can set $\beta = 1$ in the conservation laws.}
to assure that the system is scaling homogeneous
(see Section~\ref{cons-laws-scaling}).

Equations (\ref{kam-(2.18)})-(\ref{kam-(2.17)}) are a system of evolution
equations in $r$.
Noting that $\theta$ and $z$ do not explicitly appear, a {\em conservation law}
for (\ref{kam-(2.18)})-(\ref{kam-(2.17)}) simply reads
\begin{equation}
\label{conslaw}
{\mathrm D}_{r} \, T^r + {\mathrm D}_{t} \, T^t \,\dot{=} \, 0,
\end{equation}
where the total derivatives ${\mathrm D}_{r}$ and ${\mathrm D}_{t}$
are defined as
\begin{eqnarray}
\label{Dr-operator}
{\mathrm D}_{r} =  \frac{\partial}{ \partial r}
 + \sigma_r \frac{\partial} {\partial \sigma}
 + u_r \frac{\partial} {\partial u}
 + \sigma_{rr} \frac{\partial} {\partial \sigma_r}
 + u_{rr} \frac{\partial} {\partial u_r}
 + \sigma_{rt} \frac{\partial}{\partial \sigma_t}
 + u_{rt} \frac{\partial}{\partial u_t} + \ldots ,
 \end{eqnarray}
 and
\begin{eqnarray}
\label{Dt-operator}
{\mathrm D}_{t} = \frac{\partial}{ \partial t}
 + \sigma_t \frac{\partial} {\partial \sigma}
 + u_t \frac{\partial} {\partial u}
 + \sigma_{tt} \frac{\partial} {\partial \sigma_t}
 + u_{tt} \frac{\partial} {\partial u_t}
 + \sigma_{rt} \frac{\partial}{\partial \sigma_r}
 + u_{rt} \frac{\partial}{\partial u_r}
 + \ldots .
\end{eqnarray}
In (\ref{conslaw}), $T^r$ is a {\em conserved density} and $T^t$ is
the corresponding {\em flux} which are functions of $r$, $t$,
$\sigma$, and $u$, and their partial derivatives with respect to
$t$. Note that all partial derivatives of $\sigma$ and $u$ with
respect to $r$ can be eliminated by using the PDEs. The
$\,\dot{=}\,$ means that equality should only hold on solutions
$\sigma(r,t)$ and $u(r,t)$ of PDE system.

System (\ref{kam-(2.18)})-(\ref{kam-(2.17)}) has the following conservation laws
\begin{eqnarray}
\label{kam-(2.18)-(2.17)-conslaw1}
&&
{\mathrm D}_{r} (r^2 \sigma)
+ {\mathrm D}_{t} (-\delta r^2 u_t ) \,\dot{=} \, 0,
\\
\label{kam-(2.18)-(2.17)-conslaw2}
&&
{\mathrm D}_{r} (r^2 t \sigma)
+ {\mathrm D}_{t} (\delta r^2 (u - t u_t)) \,\dot{=} \, 0,
\\
\label{kam-(2.18)-(2.17)-conslaw3}
&&
{\mathrm D}_{r} (r u \sigma_t)
+ {\mathrm D}_{t}
  \Big( -\tfrac{r}{2 (n+1) \delta} \left(
  (\beta + \sigma^2)^{n+1} - (n+1) \delta^2 (u_t^2 - 2 u u_{tt}) \right)
 \Big) \,\dot{=} \, 0.
\end{eqnarray}
Conservation law (\ref{kam-(2.18)-(2.17)-conslaw1}) is (\ref{kam-(2.18)}) itself
(after multiplication by $r^2$).
Similarly, multiplying (\ref{kam-(2.18)}) by $r^2 t$ allows one to straightforwardly
recast the equation in the form (\ref{kam-(2.18)-(2.17)-conslaw2}).
Both conservation laws are easy to spot and do not require (\ref{kam-(2.17)}) to hold.
Conservation (\ref{kam-(2.18)-(2.17)-conslaw3}) can be readily verified but getting it
by inspection would require ingenuity.
It will be shown in Sections~\ref{cons-laws-scaling} and~\ref{cons-laws-multiplier}
how these conservation laws can be computed algorithmically.

Differentiating (\ref{kam-(2.18)}) with respect to $r$ once and (\ref{kam-(2.17)})
with respect to $t$ twice and equating $u_{ttr} = u_{rtt}$ yields
\begin{equation}
\label{kambapalli-wave-eq}
\sigma_{rr} + \frac{\sigma_r}{r} - \frac{4\sigma}{r^2}
= \big( \sigma (\beta + \sigma^2)^n \big)_{tt}.
\end{equation}
Conservation laws for this single wave equation for the stress will be addressed
in Section~\ref{CL-Kambapalli-Magan-general-wave-eqs}.

\subsection{Model due to Magan et al. \cite{Magan-Mason-Harley-MMS-2020} with conservation laws}
\label{model-magan-cons-laws}

The second model \cite{Magan-Mason-Harley-MMS-2020}, again in non-dimensional form,
reads
\begin{eqnarray}
\label{magan-(3.6)}
\sigma_r + \frac{\sigma}{r} &=& \delta  u_{tt},
\\
\label{magan-(3.7)}
u_r &=& \tfrac{1}{\delta} \sigma (\beta + \sigma^2)^n,
\end{eqnarray}
with the meaning of the symbols as in Section~\ref{model-kambapalli-cons-laws}.
The above system describes the propagation of axial displacement waves and shear
stress waves in a circular cylinder and a cylindrical annulus.
Based on that geometry the use of cylindrical coordinates is recommended.

Some of the conservation laws of (\ref{magan-(3.6)})-(\ref{magan-(3.7)}) are
\begin{eqnarray}
\label{magan-(3.6)-(3.7)-conslaw1}
&&
{\mathrm D}_{r} (r \sigma)
+ {\mathrm D}_{t} (-\delta r u_t ) \,\dot{=} \, 0,
\\
\label{magan-(3.6)-(3.7)-conslaw2}
&&
{\mathrm D}_{r} (r t \sigma)
+ {\mathrm D}_{t} (\delta r (u - t u_t)) \,\dot{=} \, 0,
\\
\label{magan-(3.6)-(3.7)-conslaw3}
&&
{\mathrm D}_{r} (r u \sigma_t)
+ {\mathrm D}_{t}
  \Big( -\tfrac{r}{2 (n+1) \delta} \left(
  (\beta + \sigma^2)^{n+1} - (n+1) \delta^2 (u_t^2 - 2 u u_{tt}) \right)
 \Big) \,\dot{=} \, 0.
\end{eqnarray}
Note that last conservation law coincides with (\ref{kam-(2.18)-(2.17)-conslaw3}).

Here again, by eliminating the displacement $u$, the first-order system (\ref{magan-(3.6)})-(\ref{magan-(3.7)}) can be replaced by a wave equation
for the stress
\begin{equation}
\label{magan-wave-eq}
\sigma_{rr} + \frac{\sigma_r}{r} - \frac{\sigma}{r^2}
= \big( \sigma (\beta + \sigma^2)^n \big)_{tt},
\end{equation}
for which conservation laws will be computed in Section~\ref{CL-Kambapalli-Magan-general-wave-eqs}.

\section{Computation of conservation laws using scaling homogeneity}
\label{cons-laws-scaling}

In this section we discuss the scaling-homogeneity approach
\cite{goktas-hereman-jsc-1997,hereman-etal-book-nova-2009,
hereman-etal-book-birkhauser-2005,naz-hereman-pdeam-2024,poole-hereman-jsc-2011}
for the computation of conservation laws. To keep matters
transparent we will show how to compute conservation law
(\ref{kam-(2.18)-(2.17)-conslaw3}) of
(\ref{kam-(2.18)})-(\ref{kam-(2.17)}). Since that conservation law
has an arbitrary $n$, we start with $n=1$ and repeat the
computations for $n=2$ and $3$. Once these densities and fluxes are
computed, conservation law (\ref{kam-(2.18)-(2.17)-conslaw3}) for
arbitrary $n$ is obtained by pattern matching or some
straightforward computation as shown at the end of this section.

\subsection{Scaling homogeneity}
\label{scalinghomogeneity}

One can readily verify that (\ref{kam-(2.18)})-(\ref{kam-(2.17)})
has a two-parameter family of scaling symmetries,
\begin{equation}
\label{kam-(2.18)-(2.17)-scale}
(r, t, \sigma, u, \beta) \rightarrow
(\kappa^{-(2n+1)p+q} \, r, \kappa^{-(n+1)p+q} \, t, \kappa^p \, \sigma,
\kappa^q \, u, \kappa^{2p} \, \beta),
\end{equation}
parameterized by the arbitrary real numbers $p$ and $q$.
The constant $\kappa \ne 0$ is an arbitrary scaling parameter.
If we had not introduced an auxiliary parameter $\beta$ with an
appropriate scale, then $p=0$ and (\ref{kam-(2.18)})-(\ref{kam-(2.17)})
would only have a one-parameter family of scaling symmetries.

The scaling homogeneity of (\ref{kam-(2.18)})-(\ref{kam-(2.17)}) can
be expressed in terms of weights, e.g., $W(r) = -(2n+1)p+q = -
W(D_r)$ and $\, W(t) = -(n+1)p+q = - W(D_t)$, where {\em weight}
$(W)$ of a variable (or differential operator) is the exponent of
$\kappa$ associated with that variable (operator). Thus,
(\ref{kam-(2.18)-(2.17)-scale}) is equivalent with
\begin{equation}
\label{kam-(2.18)-(2.17)-weights}
W(D_r) = (2n+1)p-q,
W(D_t) = (n+1)p-q,
W(\sigma) = p,
W(u) = q,
W(\beta) = 2p.
\end{equation}
The {\em rank} of a monomial is its total weight.
For example, $(\beta + \sigma^2)^n$ has rank $2np$.
An expression is called {\em uniform in rank} if all its
monomials have the same ranks.

If (\ref{kam-(2.18)-(2.17)-scale}) was not known it could be computed
with linear algebra as follows.
Require that (\ref{kam-(2.18)})-(\ref{kam-(2.17)}) is uniform in rank, that is,
\begin{eqnarray}
\label{kam-(2.18)-(2.17)-weighteq1}
&& W(\sigma) + W(D_r) = W(\sigma) - W(r) = W(u) + 2 W(D_t),
\\
\label{kam-(2.18)-(2.17)-weighteq2}
&& W(u) + W(D_r) = (2n+1) W(\sigma),
\quad W(\beta) = 2 W(\sigma).
\end{eqnarray}
Solve these equations for $W(D_r)$ and $W(D_t)$, to get
\begin{equation}
\label{kam-(2.18)-(2.17)-generalweights}
W(D_r) = (2n+1) W(\sigma) - W(u),
\;\;
W(D_t) = (n+1) W(\sigma) - W(u),
\end{equation}
where the arbitrary $W(\sigma) > 0 $ and $W(u) > 0 $ should be selected
so that $W(D_r)$ and $W(D_t)$ are strictly positive and preferably as small as possible.
To get the lowest possible weights, we take $W(D_r) = 1$ and $W(D_t) = n$ for
which (\ref{kam-(2.18)-(2.17)-weights}) simplifies into
\begin{equation}
\label{kam-(2.18)-(2.17)-weightssimple}
W(D_r) = n+1, \, W(D_t) = 1, \,
W(\sigma) = 1, \, W(u) = n, \, W(\beta) = 2.
\end{equation}
This simple set of weights will be used in the computations below.

\subsection{Computing conservation law (\ref{kam-(2.18)-(2.17)-conslaw3})}
\label{construct-density}

Conservation law (\ref{conslaw}) is {\em linear} in the density (and flux).
Consequently, a linear combination of densities with constant coefficients
is also a density.
Vice versa, if a density has coefficients that are, for example, powers of the parameter $\beta$, it can be split into independent densities according to those powers.
The aim is to produce linearly {\emph{independent}} densities that are as short as possible.
In particular, they should be free of constant terms and any terms that could be moved into the flux.

The key idea of the scaling-homogeneity method is that all the terms in the density
must have the same rank.
That homogeneity in rank is because (\ref{conslaw}) should only hold on solutions
of (\ref{kam-(2.18)})-(\ref{kam-(2.17)}).
Consequently, densities, fluxes, and conservation laws themselves inherit
(or adopt) the scaling homogeneity of that system (and all its other
continuous and discrete symmetries).

For example, in (\ref{kam-(2.18)-(2.17)-conslaw3}) the density has rank $1$,
the flux has rank $n+1$ and the entire conservation law has rank $n+2$.
Note that (\ref{kam-(2.18)-(2.17)-conslaw3}) is of first degree in $r$.

Using the scaling-homogeneity method and
\verb|ConservationLawsMD.m| \cite{hereman-poole-conservationlawsMD}, 
the user should select a value of $n$ (for example, $n=1$), 
specify the rank of the density to be computed
(e.g., rank one), and the desired highest degree in $r$ and $t$
(e.g., degree 1).
 \vskip 3pt \noindent {\bf Step 1}: Computation of
a candidate density $(T^r)$ for $n =1$. The code uses
(\ref{kam-(2.18)-(2.17)-weightssimple}) with $n =1$ to build a $T^r$
as a linear combination with constant coefficients of scaling
homogenous monomials involving $\beta, \sigma $, and $u$ (and their
$t$-derivatives) and the independent variables $r$ and $t$ (up to
first degree) so that each monomial has rank $1$. Table 1 shows the
terms needed together with their respective ranks.

\noindent
\begin{center}
\begin{table}[htb]
\label{kam-(2.18)-(2.17)-conslaw3-table-Tr-n1}
\begin{tabular}{ccccc}
\hline
Coefficient $r^{p_{1}} t^{p_{2}}$ & Rank & Differential term & Rank & Rank of product
\\
\hline
$1$ & $1$ & $\sigma$, $u$ & 1 & 1
\\
$r$ & $-2$ &
   $\beta\sigma$, $\sigma^3$, $\beta u$, $\sigma^2 u$, $\sigma u^2$, $u^3$, $u \sigma_t$
   & 3 & 1
\\
$t$ & $-1$ & $\sigma^2$, $\sigma u$, $u^2$ & 1 & 1
\\
\hline
\end{tabular}
\caption{
Coefficients of type $r^{p_{1}} t^{p_{2}}$ (with $p_{1}, p_{2}$ natural numbers)
of at most degree 1
(i.e., $0 \le p_{1}+p_{2} \le 1$) are paired with differential terms so that the respective products all have rank $1$.
}
\end{table}
\end{center}
\noindent The lists of differential terms are free of trivial
(constant) terms, total derivatives with respect to $t$, and
divergence-equivalent terms, meaning terms that only differ by a
$t$-derivative. For example, $\sigma u_t$ and $u \sigma_t$ of rank
$3$ are divergence-equivalent since $\sigma u_t = {\mathrm D}_{t}
(\sigma u) - u \sigma_t$. Hence, $\sigma u_t$ is removed but $u
\sigma_t$ is kept. Doing so, produces a candidate for $T^r$ which is
free of terms that could have been moved into the flux $T^t$.

Using the twelve terms in Table 1,
\begin{eqnarray}
\label{kam-(2.18)-(2.17)-conslaw3-candidate-Tr-n1}
T^r &=&
c_1 \sigma + c_2 u + c_3 \beta r \sigma + c_4 t \sigma^2
+ c_5 r \sigma^3 + c_6 \beta r u + c_7 t \sigma u + c_8 r \sigma^2 u
\\ \nonumber
&& \; + c_9 t u^2 + c_{10} r \sigma u^2 + c_{11} r u^3 + c_{12} r u \sigma_t,
\end{eqnarray}
where $c_1$ through $c_{12}$ are undetermined coefficients. See
\cite{goktas-hereman-jsc-1997,hereman-goktas-mca-2024,naz-hereman-pdeam-2024,
poole-hereman-jsc-2011} for additional details on how candidate
densities are constructed algorithmically.
\vskip 3pt \noindent {\bf
Step 2}: Computation of the undetermined coefficients. To find the
constants $c_i$, one computes
\begin{eqnarray}
\label{kam-(2.18)-(2.17)-conslaw3-DrTr-n1-org}
{\mathrm D}_r T^r
&=&
c_1 \sigma_r + c_2 u_r + c_3 \beta (\sigma + r \sigma_r)
+ 2 c_4 t \sigma \sigma_r
+ \cdots  + c_{11} u^2 (u + 3 r u_r)
\\ \nonumber
&&\, + c_{12} (u \sigma_t + r u_r \sigma_t + r u \sigma_{tr}),
\end{eqnarray}
and, using (\ref{kam-(2.18)})-(\ref{kam-(2.17)}), replaces
$\sigma_r$, $u_r$, $\sigma_{tr} = \sigma_{rt}$, to obtain
\begin{eqnarray}
\label{kam-(2.18)-(2.17)-conslaw3-DrTr-on-system-n1}
\!\!\!\!\!\!\!
P \!\!&\!=\!&\!\!
 - \tfrac{c_1}{r} ( 2 \sigma - \delta r u_{tt})
 + \tfrac{c_2}{r \delta} \left( \delta u + r \sigma (\beta + \sigma^2) \right)
 + \cdots + \tfrac{c_{12} r}{\delta}
\left(\sigma \sigma_t (\beta + \sigma^2) + \delta^2 u u_{ttt} \right)\!.
\end{eqnarray}
Since $P = {\mathrm D}_r T^r$ must match $-{\mathrm D}_t T^t$ for
some flux $T^t$ (computed below in Step 3), $P$ must be {\em exact}.
Therefore, the variational derivative (a.k.a.\ Euler operator)
\cite{hereman-etal-book-nova-2009,hereman-etal-book-birkhauser-2005}
for each of the dependent variables applied to $P$ must vanish. The
code applies the Euler operator for $\sigma$
\begin{eqnarray}
\label{euler-sigma}
{\cal{E}}_{\sigma}
&=& \sum_{k=0}^{K}
  (-{\mathrm D}_t)^k \frac{\partial }{\partial \sigma_{kt} }
  \nonumber \\
&=& \frac{\partial }{\partial \sigma}
  - {\mathrm D}_t \frac{\partial }{\partial \sigma_t}
  + {\mathrm D}_{t}^2 \frac{\partial }{\partial \sigma_{tt}}
  - {\mathrm D}_{t}^3 \frac{\partial }{\partial \sigma_{ttt}}
+ \ldots,
\end{eqnarray}
to $P$ wherein the highest derivative of $\sigma$ is $\sigma_t$.
Hence, $K =1$ and
\begin{eqnarray}
\label{eulersigma-on-P}
{\cal{E}}_{\sigma} P
\!&\!=\!&\!
\frac{\partial P}{\partial \sigma}
  - {\mathrm D}_t \frac{\partial P}{\partial \sigma_t}
\nonumber \\
\!&\!=\!&\!
-\frac{2 c_1}{r} + \frac{\beta c_2}{\delta} - \beta c_3
+ \frac{r \beta^2 c_6}{\delta} + \frac{2 t}{\delta r}(\beta r c_7- 4 \delta c_4) \sigma
+ \cdots
+ 2 \delta r c_8 u u_{tt}.
\end{eqnarray}
Next, the Euler operator for $u$ is applied to $P$ which
has a third-order term $u_{ttt}$.
Hence, $K=3$ and
\begin{eqnarray}
\label{euleru-on-P}
{\cal{E}}_{u} P
&=&
\frac{\partial P}{\partial u}
  - {\mathrm D}_t \frac{\partial P}{\partial u_t}
  + {\mathrm D}_{t}^2 \frac{\partial P}{\partial u_{tt}}
  - {\mathrm D}_{t}^3 \frac{\partial P}{\partial u_{ttt}}
\nonumber \\
&=&
\frac{c_1}{r} + 2 \beta c_6
+ \frac{t}{\delta r}(2 \beta r c_9 - \delta c_7) \sigma
+ \cdots
+ 4 \delta c_8 \sigma u_{tt} + 4 \delta r c_{10} u u_{tt}.
\end{eqnarray}
Requiring that ${\cal{E}}_{\sigma} P \equiv 0$ and ${\cal{E}}_{u} P
\equiv 0$ on the yet space (where monomials in $r, t, \sigma$, $u$,
$\sigma_t$, $u_t$, $\sigma_{tt}$, $u_{tt}$, etc., are treated as
independent) gives a {\em linear} system for the coefficients $c_i$
(not shown). Solving reveals that all $c_i$ must be zero except
$c_{12}$ which is arbitrary. Substitution of the solution and
setting $c_{12} = 1$ into
(\ref{kam-(2.18)-(2.17)-conslaw3-candidate-Tr-n1}), yields $T^r = r
u \sigma_t$ which agrees with the density in
(\ref{kam-(2.18)-(2.17)-conslaw3}). Although computed for $n=1$, in
this example the density turns out to be independent of $n$.
\vskip
3pt \noindent {\bf Step 3}: Computation of the flux $T^t$ for $n=1$.
Upon substitution of the constants $c_i$ into
(\ref{kam-(2.18)-(2.17)-conslaw3-DrTr-on-system-n1}),
\begin{equation}
\label{kam-(2.18)-(2.17)-conslaw3-DrTr-n1-eval}
P = \tfrac{r}{\delta}
\left( \sigma \sigma_t (\beta + \sigma^2) + \delta^2 u u_{ttt} \right),
\end{equation}
which must be integrated with respect to $t$ to get the flux $T^t$.
Actually, since $P = {\mathrm D}_r T^r = -{\mathrm D}_t T^t$,
after integration by parts by hand or with {\em Mathematica}, one has
\begin{equation}
\label{kam-(2.18)-(2.17)-conslaw3-Tt-n1}
T^t = -\int P \, dt
= -\tfrac{r}{4 \delta}
\left( (\beta + \sigma^2)^{2} - 2 \delta^2 (u_t^2 - 2 u u_{tt}) \right)
\end{equation}
which matches the flux in (\ref{kam-(2.18)-(2.17)-conslaw3}) for $n=1$.

In other examples, the expressions for $P$ are long and complicated.
Since repeated integration by parts with {\em Mathematica} might
fail, \verb|ConservationLawsMD.m| uses the {\em homotopy operator}
\cite[p.\ 372]{Olver-1993} to reduce the integration with respect to
$t$ to a standard integral with respect to a scaling parameter
$\lambda$. The homotopy operator turns out to be a very useful tool
for the computation of conservation laws (see, e.g.,
\cite{Anco-Bluman-EurApplMath-II-2002,hereman-etal-book-nova-2009,
hereman-etal-book-birkhauser-2005,hereman-etal-mcs-2007,hereman-goktas-mca-2024,
naz-hereman-pdeam-2024,poole-hereman-aa-2010,poole-hereman-jsc-2011}).

Application of the homotopy operator requires a few steps:
(i) computation of an integrand for each of the dependent variables
($\sigma$ and $u$),
(ii) scaling of the dependent variables (and their derivatives) with $\lambda$,
and
(iii) evaluation of one-dimensional integral with respect to $\lambda$.

In terms of the homotopy operator ${\cal H}_{{\bf u}}$,
\begin{equation}
\label{homotopysigmau-on-P}
T^t = - {\cal H}_{{\bf u}} P
    = - \int_{0}^{1} ( I_{\sigma} P + I_{u} P )[\lambda {\bf u}]
      \,\frac{d \lambda}{\lambda},
\end{equation}
where ${\bf u} = (\sigma, u)$ and $[\lambda {\bf u}]$ denotes that
(in the integrands $I_{\sigma} P$ and $I_{u} P$) $\sigma$ is
replaced by $\lambda \sigma$, $u$ by $\lambda u$, $\sigma_t$ by
$\lambda \sigma_t $, $u_t$ by $\lambda u_t$, etc. The integrand for
$\sigma$
\cite{hereman-etal-book-nova-2009,hereman-etal-mcs-2007,poole-hereman-jsc-2011}
reads
\begin{eqnarray}
\label{integrandhomotopysigma-on-P}
I_{\sigma} P
&=& \sum_{k=1}^{K}
  \left( \sum_{i=0}^{k-1} \sigma_{it} (-{\mathrm D}_t)^{k-(i+1)} \right)
  \frac{\partial P}{\partial \sigma_{kt}}
\nonumber \\
&=&
(\sigma {\mathrm I})(\frac{\partial P}{\partial \sigma_t})
+ (\sigma_t {\mathrm I}-\sigma {\mathrm D}_t)
  (\frac{\partial P}{\partial \sigma_{tt}})
+ (\sigma_{tt} {\mathrm I} - \sigma_t {\mathrm D}_t + \sigma {\mathrm D}_t^2)
     (\frac{\partial P}{\partial \sigma_{ttt}})
  + \ldots,
\end{eqnarray}
where ${\mathrm I}$ denotes the identity operator.
Applied to (\ref{kam-(2.18)-(2.17)-conslaw3-DrTr-n1-eval}) where $K=1$, one gets
\begin{equation}
\label{integrandhomotopysigma-on-P-eval}
I_{\sigma} P
= (\sigma {\mathrm I}) (\frac{\partial P}{\partial \sigma_t})
= \tfrac{r}{\delta} \sigma^2 (\beta + \sigma^2).
\end{equation}
Using (\ref{integrandhomotopysigma-on-P}) for $u$ (instead of $\sigma$) and $K=3$,
\begin{equation}
\label{integrandhomotopyu-on-P-eval} I_{u} P = - r \delta (u_t^2 - 2
u u_{tt}).
\end{equation}
Then, from (\ref{homotopysigmau-on-P}),
\begin{eqnarray}
\label{homotopysigmau-on-P-eval}
T^t &=& -\frac{r}{\delta}
\int_{0}^{1} \left( \sigma^2 (\beta + \sigma^2)
    -\delta^2 (u_t^2 - 2 u u_{tt}) \right)
  [\lambda {\bf u}] \, \frac{d \lambda}{\lambda}
  \nonumber \\
&=& - \frac{r}{\delta}
\int_{0}^{1} \lambda \left( \sigma^2 (\beta + \lambda^2 \sigma^2)
    - \delta^2 (u_t^2 - 2 u u_{tt}) \right) \, d \lambda
\nonumber \\
&=& -\tfrac{r}{4 \delta}
    \left( 2 \beta \sigma^2 + \sigma^4 - 2 \delta^2 (u_t^2 - 2 u u_{tt}) \right).
\end{eqnarray}
After adding $-\tfrac{r \beta^2}{4 \delta}$ (which does not depend
on $t$) to $T^t$, one gets (\ref{kam-(2.18)-(2.17)-conslaw3-Tt-n1}).
Once the density $T^r$ and flux $T^t$ are computed the code verifies
that they indeed satisfy (\ref{conslaw}). \vskip 3pt \noindent {\bf
Step 4}: Computation of the flux $T^t$ for arbitrary $n$. Recall
that in this example $T^r$ is valid for any value of $n$. To
determine the flux $T^t$ for arbitrary $n$, it suffices to compute
the conservation laws for $n=2$ and $n=3$ for
(\ref{kam-(2.18)})-(\ref{kam-(2.17)}) and do some pattern matching.
Or, one can provide the form of $T^r$ (e.g., $c_1 r u \sigma_t$) and
let \verb|ConservationLawsMD.m| compute $T^t$ automatically.

Alternatively, one can evaluate
\begin{equation}
\label{kam-(2.18)-(2.17)-conslaw3-DrTr-n}
{\mathrm D}_r T^r
= {\mathrm D}_r (r u \sigma_t)
= u \sigma_t + r u_r \sigma_t + r u \sigma_{tr}
\end{equation}
on solutions of (\ref{kam-(2.18)})-(\ref{kam-(2.17)}) yielding
\begin{equation}
\label{kam-(2.18)-(2.17)-conslaw3-DrTr-on-system-n}
P = \tfrac{r}{\delta} \left(\sigma \sigma_t
(\beta + \sigma^2)^n + \delta^2 u u_{ttt} \right).
\end{equation}
Thus,
\begin{equation}
\label{kam-(2.18)-(2.17)-conslaw3-Tt-n}
T^t = -\int P \, dt
= -\tfrac{r}{2 (n+1) \delta}
  \left( (\beta + \sigma^2)^{n+1} - 2 (n+1) \delta^2 (u_t^2 - 2 u u_{tt}) \right)
\end{equation}
which matches the flux in (\ref{kam-(2.18)-(2.17)-conslaw3}).

Conducting a search with \verb|ConservationLawsMD.m| for
conservation laws of (\ref{kam-(2.18)})-(\ref{kam-(2.17)}) and
(\ref{magan-(3.6)})-(\ref{magan-(3.7)}), respectively, where the
densities have coefficients up to {\emph{third degree}} in $r$ and
$t$, yields
(\ref{kam-(2.18)-(2.17)-conslaw1})-(\ref{kam-(2.18)-(2.17)-conslaw3})
and
(\ref{magan-(3.6)-(3.7)-conslaw1})-(\ref{magan-(3.6)-(3.7)-conslaw3}),
respectively. No additional conservation laws for either model could
be found. This is in stark contrast with the successful computation
of seven conservation laws for a similar model in Cartesian
coordinates \cite{naz-hereman-pdeam-2024} which is conjectured to
have infinitely many conservation laws.

\section{$\!\!$Computation of conservation laws with a multiplier method}
\label{cons-laws-multiplier}

Several methods are available to compute conservation laws (see,
e.g.,
\cite{Bluman-Anco-book-2002,Bluman-Cheviakov-Anco-book-2010,Naz-Mahomed-Mason-AMC-2008}).
We refer the reader to
\cite{Anco-FIC-chapter-2017,Cheviakov-JEM-2010,Naz-Mahomed-Mason-AMC-2008,Wolf-EJAM-2002}
for reviews and comparisons of commonly-used approaches.

In this section, we discuss the multiplier method
\cite{Bluman-Anco-book-2002,Naz-Mahomed-Mason-AMC-2008,Olver-1993,Steudel-1962}
for the computation of conservation laws. Using the multiplier
approach, we derive the conservation laws
(\ref{kam-(2.18)-(2.17)-conslaw1})-(\ref{kam-(2.18)-(2.17)-conslaw3})
corresponding to the system (\ref{kam-(2.18)})-(\ref{kam-(2.17)}).
Similarly, the conservation laws
(\ref{magan-(3.6)-(3.7)-conslaw1})-(\ref{magan-(3.6)-(3.7)-conslaw3})
for the system (\ref{magan-(3.6)})-(\ref{magan-(3.7)}) are also
obtained. Instead of treating these systems separately, we combine
them in a single system of PDEs by introducing suitable constants.
Furthermore, by considering a general function $F(\sigma)$, this
set-up covers a broad class of models for wave propagation in
elastic materials formulated in cylindrical coordinates.

\subsection{Conservation laws for (\ref{kam-(2.18)})-(\ref{kam-(2.17)}) and (\ref{magan-(3.6)})-(\ref{magan-(3.7)})}

Model (\ref{kam-(2.18)})-(\ref{kam-(2.17)}) by Kambapalli et al.\
\cite{Kambapalli-Kannan-Rajagopal-2014} and system
(\ref{magan-(3.6)})-(\ref{magan-(3.7)}) by Magan et al.\
\cite{Magan-Mason-Harley-MMS-2020} can be combined and generalized
as follows:
\begin{eqnarray}
\label{general1}
\sigma_r + \frac{\kappa_1 \sigma}{r} &=& \delta  u_{tt},
\\
\label{general2}
u_r + \frac{\kappa_2 u}{r} &=& \tfrac{1}{\delta} F(\sigma).
\end{eqnarray}
Eliminating $u$ from (\ref{general1})-(\ref{general2}) yields
\begin{equation}
\label{Kambapalli_Magan_Combined_PDE}
 \sigma_{rr} + (\kappa_1+\kappa_2) \frac{\sigma_r}{r}
   - \kappa_1 (1-\kappa_2) \frac{\sigma}{r^2}
 = F(\sigma)_{tt}.
\end{equation}
For $F(\sigma) = \sigma (\beta + \sigma^2)^n$ and
$\kappa_1 = 2$ and $\kappa_2 = -1$,
system (\ref{general1})-(\ref{general2}) matches (\ref{kam-(2.18)})-(\ref{kam-(2.17)})
and (\ref{Kambapalli_Magan_Combined_PDE}) becomes (\ref{kambapalli-wave-eq}).
Likewise, for $\kappa_1 = 1$ and $\kappa_2 = 0$, one gets (\ref{magan-(3.6)})-(\ref{magan-(3.7)}) and (\ref{magan-wave-eq}).
Note that $\kappa_1 + \kappa_2 = 1$ for either set of parameters.

Several constitutive relations, $\epsilon = F(\sigma)$, relating
scalar linearized strain ($\epsilon$) to Cauchy stress ($\sigma$,
also scalar) have been reported in the literature. A dozen of those
are summarized in Table 2 together with references where additional
information can be obtained.

The examples are meant to illustrate the form $F(\sigma)$ can take.
Not all these functions $F(\sigma)$ might be relevant for
(\ref{general1})-(\ref{general2}). One should check the cited
references to derive the appropriate one-dimensional version of the
implicit constitutive relation,
${\boldsymbol{\epsilon}} = {\bf F}({\boldsymbol{\sigma}})$,
relating the linearized strain
(${\boldsymbol{\epsilon}}$) and Cauchy stress
(${\boldsymbol{\sigma}}$) tensors for purely elastic deformations.
It is also important to consider the most suitable coordinate system
for the geometrical configuration and derive the correct equations
of motion for the problem at hand. The reader is referred to
\cite{bustamante-in-book-2026,bustamante-rajagopal-in-book-2020} for
a discussion about the appropriateness and limitations of the
constitutive relations within their physical contexts.

\vfill
\newpage
\begin{center}
\begin{table}[ht!]
\label{F-combined-models}
\small
\begin{tabular}{lll}
\hline
$F(\sigma)$ & $\!\!\!\!\!\!\!$ Remarks & $\!\!\!$ References
\\
\hline
\\
$ \sigma \big( 1 + \sigma^2 \big)^n $
& $\!\!\!\!\!\!\!$ $n \ge 0$ rational
& $\!\!\!$
\cite[[Eq.\ (2.17)]{Kambapalli-Kannan-Rajagopal-2014},
\cite[Eq.\ (28)]{kannan-rajagopal-saccomandi-wm-2014},
\\
& & $\!\!\!$ \cite[Eqs.\ (4.2.7)\;\&\;(5.2.6)]{magan-thesis-2018},
\\
& & $\!\!\!$
\cite[Eq.\ (2.20)]{magan-etal-wm-2018},
\cite[Eq.\ (1)]{naz-hereman-pdeam-2024}
\\
\\
$ \frac{1}{\delta}
  \Big( \beta \sigma + \alpha \sigma (1 + \frac{\gamma}{2} \sigma^2)^n \Big) $
& $\!\!\!\!\!\!\!$ $n$ rational
& $\!\!\!$
\cite[Eq.\ (2.4)]{Erbay-Sengul-2015},
\cite[Eq.\ (28)]{kannan-rajagopal-saccomandi-wm-2014}
\\
\\
$ \frac{\alpha \delta}{2}
  \sigma \big(1 + \frac{1}{1 + \sigma^2} \big)^n $
& $\!\!\!\!\!\!\!$ $n \ge 0$ rational & $\!\!\!$ \cite[Eqs.\
(4.7.4)\;\&\;(5.6.2)]{magan-thesis-2018}, \cite[Eq.\
(2.21)]{magan-etal-wm-2018},
\\
& & $\!\!\!$
\cite[Eq.\ (2.9)]{Magan-Mason-Harley-MMS-2020},
\cite[Eq.\ (2.17)]{rajagopal-MMS-2011}
\\
\\
$ \frac{\alpha \delta}{2} \sigma (\frac{1}{1 + \sigma^2})^n $
& $\!\!\!\!\!\!\!$ $n \ge \frac{1}{2}$, rational
& $\!\!\!$
\cite[Eq.\ (2.10)]{Magan-Mason-Harley-MMS-2020}
\\
\\
$ \frac{\alpha}{\delta} \sigma
\left( \frac{1}{(1 + |\sigma|^n)^{\frac{1}{n}}} \right) $
& $\!\!\!\!\!\!\!$ $n > 0$, integer
& $\!\!\!$
\cite[Eq.\ (2.5)]{Erbay-Sengul-2015}
\cite[Eq.\ (3.10), $n = 1$]{rajagopal-MMS-2011}
\\
\\
$ \frac{\alpha}{\delta} \Big( -1 + \frac{1}{1 + \beta \sigma}
  + \frac{\gamma \sigma}{\sqrt{1 + \iota \sigma^2}} \Big) $
& & $\!\!\!$
\cite[Eq.\ (25)]{bustamante-rajagopal-ijnlm-2011},
\cite[Eq.\ (29)]{bustamante-sfyris-MMS-2015},
\\
& & $\!\!\!$
\cite[Eq.\ (15)]{ortiz-bernardin-etal-IJSS-2014},
\cite[Eq.\ (3.2)]{Sengul-AES-2021}
\\
\\
$ \frac{\alpha}{\delta} \Big(- \frac{\beta \sigma}{1 + \beta \sigma}
  + \frac{\mu \sigma}{\sqrt{1 + \gamma^2 \sigma^2}} \Big) $
& & $\!\!\!$
\cite[Eq.\ (47)]{bridges-rajagopal-ZAMP-2015}
\\
\\
$ \frac{\alpha}{\delta}
\Big(1-\frac{1}{1+\frac{\sigma}{1+\tilde{\delta} |\sigma|}} \Big) $
&  &
\\
$ + \frac{\beta}{\delta} \sigma \big(1 + \frac{1}{1 + \gamma \sigma^2} \big)^n $
& $\!\!\!\!\!\!\!$ $n$ rational
& $\!\!\!$
\cite[Eq.\ (2.7)]{Erbay-Sengul-2015},
\cite[Eq.\ (3.13)]{rajagopal-MCA-2010}
\\
\\
$ \frac{1}{\delta} {\mathcal{N}} {\mathrm{arctan}} (\theta \sigma)$
& & $\!\!\!$
\cite[Eq.\ (3.1)]{Sengul-AES-2021}
\\
\\
$ \frac{1}{\delta}
  \big(\!-\alpha \,{\mathrm{tanh}} (\beta \sigma)
  + \gamma \sigma \frac{1}{\sqrt{1 + \iota \sigma^2}} \!\big) $
& & $\!\!\!$
\cite[Eq.\ (13)]{bustamante-etal-imajam-2024}
\\
\\
$ \frac{\alpha}{\delta} \Big( 1 -{\mathrm{e}}^{-\beta \sigma}
  + \frac{\gamma \sigma}{1 + \iota |\sigma|} \Big) $
& & $\!\!\!$
\cite[Eq.\ (28)]{bustamante-rajagopal-ijnlm-2011}
\\
\\
$ \frac{1}{\delta} \Big( \beta_2 \sigma
+ \beta_3 \sigma {\mathrm{e}}^{(1 + \beta_4 \sigma^2)^{\tfrac{n}{2}}} \Big)$
& $\!\!\!\!\!\!\!$ $n\ge 0$, integer &  $\!\!\!$
\cite[Eq.\ (28)]{devendiran-etal-IJSS-2017}
\\
\\
$ \frac{\alpha \gamma}{\delta}
  \big( \frac{\sigma}{1 + |\sigma|} \big)$
& &
\\
$ + \frac{\alpha}{\delta}
  \Big(1-{\mathrm{e}}^{-\frac{\beta \sigma}{1 + \tilde{\delta} |\sigma|}} \Big)$
& & $\!\!\!$
\cite[Eq.\ (2.6)]{Erbay-Sengul-2015},
\cite[Eq.\ (3.12)]{rajagopal-MCA-2010}
\\
\\
$ \frac{\beta}{\delta} \sigma
  \big( (1 + \gamma \sigma^n)^{-\frac{1}{n}} \big) $
& &
\\
$ + \frac{\alpha}{\delta}
 \Big( 1-{\mathrm{e}}^{-\frac{\lambda \sigma}{1 + \tilde{\delta}\sigma}}\Big)$
& $\!\!\!\!\!\!\!$ $n$ rational
& $\!\!\!$
\cite[Eq.\ (2.6), n = 1]{Erbay-Sengul-2015},
\cite[Eq.\ (28)]{Kulvait-etal-IJF-2013},
\\
& &
$\!\!\!$
\cite[Eq.\ (3.14)]{rajagopal-MCA-2010},
\cite[Eq.\ (34), n = 1]{rajagopal-MMS-2011-conspectus-paper}
\\
\\
\hline
\end{tabular}
\caption{
Examples of scalar (1D) versions of constitutive relations from the literature.
\newline
The meaning of the material parameters ($\alpha, \beta, \gamma,
\iota $, etc.) can be found in the cited references. }
\end{table}
\end{center}
\vfill
\newpage

Therefore we investigate system (\ref{general1})-(\ref{general2})
for a class of constitutive relations with an unspecified
$F(\sigma)$ and arbitrary $\kappa_1$ and $\kappa_2$. System
(\ref{general1})-(\ref{general2}) can be expressed as
\begin{equation}
\label{CL-Kambapalli-model}
E^1(r,t,\sigma,u,\sigma_{r},u_{tt}) = 0,
\;\;
E^2(r,t,\sigma,u,u_{r}) = 0,
\end{equation}
where
\begin{equation}
\label{gen-model}
E^1 = \sigma_r+\frac{\kappa_1 \sigma}{r} -\delta u_{tt},
\;\;
E^2 =u_r+\frac{\kappa_2 u}{r} - \tfrac{1}{\delta} F(\sigma).
\end{equation}
The conservation law is written in characteristic form
\cite{Naz-Mahomed-Mason-AMC-2008,Olver-1993,Steudel-1962} as
\begin{equation}
\label{characteristics_equation_gen}
{\mathrm D}_{r} T^r + {\mathrm D}_{t} T^t
= \Lambda^1 \left(\sigma_r+\frac{\kappa_1 \sigma}{r}- \delta u_{tt}\right)
  + \Lambda^2 \left(u_r+\frac{\kappa_2 u}{r} -\tfrac{1}{\delta} F(\sigma)\right),
\end{equation}
with ${\mathrm D}_{r}$ and ${\mathrm D}_{t}$ given in
(\ref{Dr-operator}) and (\ref{Dt-operator}). The differential
functions $\Lambda^1$ and $\Lambda^2$ are the multipliers (a.k.a.\
characteristics). We assume $\Lambda^1(r,t,\sigma, u,\sigma_t, u_t)$
and the same dependencies for $\Lambda^2$. Once restricted to
{\emph{solutions}} of (\ref{CL-Kambapalli-model}), the
characteristic equation (\ref{characteristics_equation_gen}) yields
a local conservation law (\ref{conslaw}). As outlined in Olver
\cite{Olver-1993}, the determining equations for the multipliers
$\Lambda^1$ and $\Lambda^2$ are established by applying the Euler
operators to the characteristic equation
(\ref{characteristics_equation_gen}) with respect to the dependent
variables $\sigma$ and $u$. Hence,
\begin{equation}
\label{multipliers_eq1_gen}
\mathcal{E}_{\sigma} \left[ \Lambda^1 \left(\sigma_r+\frac{\kappa_1 \sigma}{r}
- \delta u_{tt}\right)+\Lambda^2 \left(u_r+\frac{\kappa_2 u}{r}
- \tfrac{1}{\delta} F(\sigma)\right) \right] = 0,
\end{equation}
\begin{equation}
\label{multipliers_eq2_gen}
\mathcal{E}_u\left[ \Lambda^1 \left(\sigma_r+\frac{\kappa_1 \sigma}{r}
- \delta u_{tt}\right)+\Lambda^2 \left(u_r+\frac{\kappa_2 u}{r}
- \tfrac{1}{\delta} F(\sigma)\right)\right] = 0,
\end{equation}
where $\mathcal{E}_{\sigma}$ and $\mathcal{E}_u$ are the standard
Euler operators 
(see (\ref{euler-sigma}) and (\ref{euleru-on-P}), respectively). 
Equations (\ref{multipliers_eq1_gen}) and (\ref{multipliers_eq2_gen}) 
result in an overdetermined system of PDEs for the multipliers 
$\Lambda^1$ and $\Lambda^2$. 
The derivation and solution of such determining equations has been 
implemented in \verb|GeM|
\cite{Cheviakov-CPC-2007,Cheviakov-JEM-2010,Cheviakov-website-2025} 
and other symbolic packages \cite{Wolf-EJAM-2002}.

After some simplifications, the determining equations are
\begin{eqnarray}
\label{overdetermined_system_multipliers_general}
\Lambda^1_{\sigma_t} = 0, \;
\Lambda^1_{u_t u_t} = 0, \;
\Lambda^1_{\sigma} = 0, \;
\Lambda^2_{u_t} = 0, \;
\Lambda^2_{\sigma_t\sigma_t} = 0,
\nonumber \\
\Lambda^1_{u_t}+\Lambda^2_{\sigma_t} = 0, \;
\Lambda^1_u-\Lambda^2_{\sigma}-\Lambda^1_{tu_t}-u_t \Lambda^1_{u u_t} = 0,
\nonumber \\
\Lambda^1_u-\Lambda^2_{\sigma}+\Lambda^2_{t\sigma_t}
+ u_t\Lambda^2_{u\sigma_t}+\sigma_t\Lambda^2_{\sigma\sigma_t} = 0,
\nonumber \\
2\Lambda^1_u+\Lambda^1_{tu_t}+u_t \Lambda^1_{uu_t} = 0,
\nonumber \\
\left(\Lambda^{2}_{t\sigma_t} -\Lambda^{2}_{\sigma}
+ u_t \Lambda^{2}_{u \sigma_t}
+ \sigma_t \Lambda^{2}_{\sigma\sigma_t}\right)F(\sigma)
+ \left(\sigma_t \Lambda^{2}_{\sigma_t}-\Lambda^{2}\right) F'(\sigma)
\\
+\tfrac{\delta}{r} \left[k_2u \left(\Lambda^{2}_{\sigma}
 -u_t \Lambda^{2}_{u\sigma_t} -\sigma_t \Lambda^{2}_{\sigma \sigma_t}
 -\Lambda^{2}_{t\sigma_t}\right)+k_1 \Lambda^{1}
 - r \Lambda^{1}_{r}-k_2 u_t \Lambda^{2}_{\sigma_t} \right] = 0,
\nonumber\\
\Lambda^{2}_{u}F(\sigma)
 + \delta \Lambda^{2}_{r}
 + \delta^2 ( \Lambda^{1}_{uu} u_t^{2}
 + \Lambda^{1}_{tt}+2 \Lambda^{1}_{t u} u_t )
 \nonumber\\
 -\tfrac{\delta}{r} \left[ k_1\sigma\left( \Lambda^{1}_{u}
 - u_t \Lambda^{1}_{u_tu} -\Lambda^{1}_{u_tt} \right)
 + k_2 u \Lambda^{2}_{u}+k_2 \Lambda^2
 - k_1\sigma_t \Lambda^{1}_{u_t} \right] = 0.
\nonumber
\end{eqnarray}
Solving (\ref{overdetermined_system_multipliers_general}) for an
arbitrary $F(\sigma)$ yields the following multipliers
\begin{equation}
\label{gen_multipliers}
\Lambda^1 = c_1 r^{\kappa_1} + c_2 t r^{\kappa_1}
            - c_3 r^{\kappa_1 + \kappa_2} u_t,
\quad
\Lambda^2 = c_3 r^{\kappa_1 + \kappa_2} \sigma_t.
\end{equation}
Next, one substitutes (\ref{gen_multipliers}) into
(\ref{characteristics_equation_gen}) and determines $T^r$ and $T^t$.
Direct integration \cite{Cheviakov-Zhao-book-2024} as implemented in \verb|GeM|
\cite{Cheviakov-CPC-2007,Cheviakov-JEM-2010,Cheviakov-website-2025} 
yields 
\begin{eqnarray}
\label{CL_gen}
&& {\mathrm D}_{r} \left[ c_1  r^{\kappa_1} \sigma
 + c_2 r^{\kappa_1}  t \sigma+c_3 r^{\kappa_1+\kappa_2} u \sigma_t \right]
 + {\mathrm D}_{t} \left[ - c_1 \delta r^{\kappa_1} u_t
 + c_2 \delta r^{\kappa_1} (u- t u_t) \right.
\nonumber \\
&&
\left. + c_3\, r^{\kappa_1+\kappa_2} \bigg( \tfrac{\delta}{2}(u_t^2-2uu_{tt})
 - \tfrac{1}{\delta} \int F(\sigma)\, d\sigma \bigg) \right]
 \,\dot{=} \, 0.
\end{eqnarray}
Alternatively, one could use a homotopy operator method
\cite{Anco-Bluman-EurApplMath-II-2002,poole-hereman-jsc-2011} to
compute $T^r$ and $T^t$.

Any density-flux pair for system (\ref{general1})-(\ref{general2})
with multipliers of the form $\Lambda^1(r,t,\sigma, u, \sigma_t,
u_t)$ and  $\Lambda^2(r,t,\sigma, u, \sigma_t, u_t)$ is therefore a
linear combination of the three conserved densities and associated
fluxes presented in Table 3 together with the corresponding
multipliers. For the function $F(\sigma)$ one could take examples
from Table 2.

\begin{center}
\begin{table}[h!]
\label{Kambapalli_Magan_Combined_SYS_CLs}
\begin{tabular}{lll}
\hline Multipliers ($\Lambda_1, \Lambda_2$) &\quad Density $T^r$ & \quad Flux $T^t$
\\ \hline
\\
$(r^{\kappa_1}$, $0$) & \quad $r^{\kappa_1} \sigma$ & \quad
$-\delta
r^{\kappa_1} u_t$
\\
\\
$ ( t r^{\kappa_1}$, $0$) & \quad $ r^{\kappa_1} t \sigma $ & \quad
$ \delta r^{\kappa_1} (u -t u_t)$
\\
\\
$ ( - r^{\kappa_1 + \kappa_2} u_t$, $ r^{\kappa_1 + \kappa_2}
\sigma_t $) & \quad $ r^{\kappa_1+\kappa_2} u \sigma_t$ & \quad
$
r^{\kappa_1+\kappa_2} \big(\tfrac{\delta}{2}(u_t^2 -2 u u_{tt})
       - \tfrac{1}{\delta} \int F(\sigma)\, d\sigma \big)$
\\
\\
\hline
\\
\end{tabular}
\caption{ Multipliers, conserved densities, and fluxes for system
(\ref{general1})-(\ref{general2}) where $\kappa_1$ and $\kappa_2$
are arbitrary. The densities for
(\ref{kam-(2.18)})-(\ref{kam-(2.17)}) corresponding to $\kappa_1 =
2$ and $\kappa_2 = -1$ are given in
(\ref{kam-(2.18)-(2.17)-conslaw1})-(\ref{kam-(2.18)-(2.17)-conslaw3}).
Those for (\ref{magan-(3.6)})-(\ref{magan-(3.7)}) corresponding to
$\kappa_1 = 1$ and $\kappa_2 = 0$ can be found in
(\ref{magan-(3.6)-(3.7)-conslaw1})-(\ref{magan-(3.6)-(3.7)-conslaw3}).
}
\end{table}
\end{center}

In summary, once the dependencies of the multipliers are selected,
the code \verb|GeM| computes the multipliers (\ref{gen_multipliers})
as well as the densities and associated fluxes using one of several
approaches \cite{Cheviakov-JEM-2010}: direct integration
\cite{Cheviakov-Zhao-book-2024}, homotopy techniques
\cite{Anco-Bluman-EurApplMath-II-2002,hereman-etal-book-birkhauser-2005,
hereman-etal-mcs-2007,poole-hereman-jsc-2011} or a scaling-symmetry
based formula for the conserved vector \cite{anco-jpa-2003}.
Instructions for how to use \verb|GeM| can be found on Cheviakov's
website \cite{Cheviakov-website-2025}.

Restricted to the models in Section~\ref{model-eqs-cons-laws}  for
which $\kappa_1 + \kappa_2 = 1$, with \verb|ConservationLawsMD.m| we
computed conservation laws of  (\ref{general1})-(\ref{general2})
with $F(\sigma) = \sigma (\beta + \sigma^2)^n$ and $\kappa_2 =
1-\kappa_1$. Doing so for $n = 1$ and coefficients up to degree
three in $r$ and $t$, the output of \verb|ConservationLawsMD.m|
confirms the results in Table 3.

\subsection{Conservation laws of (\ref{Kambapalli_Magan_Combined_PDE})} \label{CL-Kambapalli-Magan-general-wave-eqs}

For completeness, we compute low-order conservation laws for the
wave equation (\ref{Kambapalli_Magan_Combined_PDE}) associated with
(\ref{general1})-(\ref{general2}). The results for
(\ref{kambapalli-wave-eq}) and (\ref{magan-wave-eq}) associated with
(\ref{kam-(2.18)})-(\ref{kam-(2.17)}) and
(\ref{magan-(3.6)})-(\ref{magan-(3.7)}), respectively, are special
cases. The results are summarized in Table 4.

\begin{center}
\begin{table}[h]
\label{Table-Kambapalli_Magan_Combined_PDE_CLs}
\begin{tabular}{llll}
\hline
Multiplier $\Lambda$ & Density $T^r$ & \; Flux $T^t$
\\
\hline
\\
$r^{\kappa_2}$ & $ r^{\kappa_2 -1}(r \sigma_r + \kappa_1 \sigma)$ &
\;
$ -r^{\kappa_2} \sigma_t F^{\prime}$
\\
\\
$t r^{\kappa_2}$ & $ t r^{\kappa_2 -1}(r \sigma_r + \kappa_1
\sigma)$ & \;
$ r^{\kappa_2}(F-t\sigma_t F^{\prime})$
\\
\\
$r^{1+\kappa_1}$ & $ r^{\kappa_1} \big(r \sigma_r - (1-\kappa_2)
\sigma \big)$ & \;
$ -r^{1+\kappa_1} \sigma_t F^{\prime}$
\\
\\
$ t r^{1+\kappa_1}$ & $ t r^{\kappa_1} \big( r\sigma_r -(1-\kappa_2)
\sigma \big)$ & \;
$r^{\kappa_1 + 1}(F - t\sigma_t F^{\prime})$
\\
\\
\hline
\\
$ \frac{1}{r}$ & $ \frac{1}{r^2}(r\sigma_r+2 \sigma)$ & \;
$
-\frac{1}{r} \sigma_t (\beta + \sigma^2)^{n-1}
      \big( \beta + (2n+1)\sigma^2 \big) $
\\
\\
$ \frac{t}{ r}$ & $ \frac{t}{r^2}( r\sigma_r+2 \sigma)$ & \;
$
\frac{1}{r} (\beta+\sigma^2)^{n-1} \big(\sigma (\beta+\sigma^2)
     -t \sigma_t \big( \beta + (2n+1)\sigma^2 \big)\big)$
\\
\\
$r^3$ & $ r^2(r\sigma_r-2 \sigma)$ & \;
$ -r^3 \sigma_t (\beta
+\sigma^2)^{n-1} \big(\beta +(2n+1)\sigma^2\big)$
\\
\\
$t r^3$ & $t r^2 ( r\sigma_r-2 \sigma) $ & \;
$r^3
(\beta+\sigma^2)^{n-1} \big( \sigma (\beta+\sigma^2)
    -t \sigma_t \big( \beta+(2n+1)\sigma^2 \big) \big)$
\\
\\
\hline
\\
$1$ & $ \frac{1}{r}(r\sigma_r+\sigma)$ & \;
$ -\sigma_t (\beta
+\sigma^2)^{n-1} \big(\beta +(2n+1)\sigma^2\big)$
\\
\\
$ t $ & $ \frac{t }{r}( r\sigma_r+\sigma)$ & \;
$(\beta+\sigma^2)^{n-1} \big( \sigma(\beta+\sigma^2)
  - t \sigma_t \big(\beta+(2n+1)\sigma^2 \big) \big)$
\\
\\
$r^2$ & $ r(r\sigma_r- \sigma)$ & \;
$ -r^2\sigma_t (\beta
+\sigma^2)^{n-1} \big(\beta + (2n+1)\sigma^2\big)$
\\
\\
$t r^2$ & $ t r(r\sigma_r- \sigma)$ & \;
$ r^2
(\beta+\sigma^2)^{n-1} \big( \sigma(\beta+\sigma^2)
    -t \sigma_t \big(\beta+(2n+1)\sigma^2 \big) \big)$
\\
\\
\hline
\end{tabular}
\caption{ The top set are the multipliers, conserved densities, and
fluxes for the parameterized wave equation
(\ref{Kambapalli_Magan_Combined_PDE}) with arbitrary $\kappa_1$ and
$\kappa_2$. The middle set is for (\ref{kambapalli-wave-eq}) where
$\kappa_1 = 2$ and $\kappa_2 = -1$. The bottom set is for
(\ref{magan-wave-eq}) where $\kappa_1 = 1$ and $\kappa_2 = 0$. }
\end{table}
\end{center}

\section{Application 1: Khokhlov-Zabolotskaya-Kuznetsov equation}
\label{application-KZK} In this section we investigate the
Khokhlov-Zabolotskaya-Kuznetsov (KZK) equation which is primarily
used to model the propagation of sound beams in various nonlinear
media. We compute conservation laws for the KZK equation in
Cartesian coordinates as well as cylindrical and spherical
coordinates.

\subsection{Conservation laws of the Khokhlov-Zabolotskaya-Kuznetsov equation
in Cartesian coordinates}
\label{KZK-equation-cartesian}

The KZK equation
\cite[Ch.\ 3, Eq.\ (3.65)]{Hamilton-Morfey-ch3-springer-2024},
\begin{equation}
\label{KZK-org}
p_{zt} - \tfrac{c_0}{2} \nabla_{\perp}^2 p - \tfrac{\delta}{2 c_0^3} p_{ttt}
- \tfrac{\beta}{2 \rho_0 c_0^3} (p^{2})_{tt} = 0,
\end{equation}
describes the propagation of a beam in a nonlinear medium. For
simplicity of notation we have replaced the retarded time variable
$(\tau)$ by $t$. The second term in (\ref{KZK-org}) accounts for the
diffraction of the beam. The operator $\nabla_{\perp}^2 =
\frac{\partial^2}{\partial x^2}+\frac{\partial^2}{\partial y^2}$ is
the Laplacian acting on the plane perpendicular to the $z-$axis
along which the beam propagates. For acoustic waves, $p$ denotes the
sound pressure, $\rho_0$ a reference density of the medium, and
$c_0$ the speed of sound. The positive parameters $\beta$ and
$\delta$ are coefficients of nonlinearity and diffusion (viscosity,
absorption), respectively. The reader is referred to
\cite{Rozanova-Pierrat-CMS-2009} for a detailed discussion of
various derivations of (\ref{KZK-org}) and historical references.

Clearing denominators yields
\begin{equation}
\label{KZK}
\delta p_{ttt} + {\tilde{\beta}} (p^{2})_{tt}
 - 2 c_0^3 p_{zt} + c_{0}^4 (p_{xx} + p_{yy}) = 0,
\end{equation}
where $\tilde{\beta} = {\frac{\beta}{\rho_0}}$, for which we will
compute low-order conservation laws with the multiplier method. The
conservation law in characteristic form reads as
\begin{equation}
\label{characteristics_equation_KZK}
{\mathrm D}_{t} \, T^t + {\mathrm D}_{x} \, T^x
+ {\mathrm D}_{y} \, T^y + {\mathrm D}_{z} \,T^z
= \Lambda \left( \delta p_{ttt} + {\tilde{\beta}} (p^{2})_{tt}
  -2 c_0^3 p_{zt} + c_{0}^4 (p_{xx}+p_{yy}) \right),
\end{equation}
where $\Lambda$ is a characteristic (a.k.a.\ multiplier). Total
derivative ${\mathrm D}_{t}$ is given in (\ref{Dt-operator}) and
${\mathrm D}_{x},\; {\mathrm D}_{y}$, and ${\mathrm D}_{z}$ are
defined analogously. When (\ref{characteristics_equation_KZK}) is
evaluated on solutions of (\ref{KZK}) one gets a local conservation
law. To compute low-order conservation laws we assume
$\Lambda(x,y,z,t,p,p_x,p_y,p_z,p_t)$.

Taking the variational derivative of the characteristic equation
(\ref{characteristics_equation_KZK}) yields
\begin{equation}
\label{multipliers_eq_KZK}
{\mathcal{E}_{p}}
\left[ \Lambda \left(\delta p_{ttt}
 + {\tilde{\beta}} (p^{2})_{tt} -2 c_0^3 p_{zt}
 + c_{0}^4 (p_{xx}+p_{yy}) \right)
\right] = 0,
\end{equation}
where $\mathcal{E}_{p}$ is the Euler operator
\cite{poole-hereman-jsc-2011}
\begin{eqnarray}
\label{Euler-operator-4D}
{\mathcal{E}_{p}}
&=& \sum_{k=0}^{K} \sum_{\ell=0}^{L} \sum_{m=0}^{M} \sum_{n=0}^{N}
(-{\mathrm D}_x)^{k}
(-{\mathrm D}_y)^{\ell}
(-{\mathrm D}_z)^{m}
(-{\mathrm D}_t)^{n}
\frac{\partial }{\partial p_{k x \, \ell y \, m z \, n t}}
\nonumber \\
&=&
  \frac{\partial }{\partial p}
  - {\mathrm D}_x \frac{\partial }{\partial p_x}
  - {\mathrm D}_y \frac{\partial }{\partial p_y}
  - {\mathrm D}_z \frac{\partial }{\partial p_z}
  - {\mathrm D}_t \frac{\partial }{\partial p_t}
  +  {\mathrm D}_{x}^2 \frac{\partial }{\partial p_{xx}}
  + {\mathrm D}_{y}^2 \frac{\partial }{\partial p_{yy}}
\nonumber \\
&& + \, {\mathrm D}_{z}^2 \frac{\partial }{\partial p_{zz}}
   +  {\mathrm D}_{t}^2 \frac{\partial }{\partial p_{tt}}
   + \cdots
   + {\mathrm D}_{z}\, {\mathrm D}_{t} \frac{\partial }{\partial p_{zt}}
   + \cdots
   - {\mathrm D}_{t}^3 \frac{\partial }{\partial p_{ttt}}
   + \ldots,
\end{eqnarray}
where $K$, $L$, $M$ and $N$ are the highest-orders of respective
derivatives of $p$ appearing in the expression ${\mathcal{E}_{p}}$
acts upon.

Equation (\ref{multipliers_eq_KZK}) results in an
overdetermined system of PDEs for the multiplier $\Lambda$ (which
can be automatically computed with \verb|GeM|):
\begin{equation}
\label{overdetermined_system_multiplier_KZK}
 \Lambda_{tt} = 0, \;\,
 \Lambda_{p} = 0, \;\,
 \Lambda_{p_x} = 0, \;\,
 \Lambda_{p_y} = 0, \;\,
 \Lambda_{p_z} = 0, \;\,
 \Lambda_{p_t} = 0, \;\,
 \Lambda_{tz} = \tfrac{c_0}{2} \left(\Lambda_{xx} + \Lambda_{yy}
 \right).
\end{equation}
Solving (\ref{overdetermined_system_multiplier_KZK}) yields
\begin{equation}
\label{gen_multipliers_KZK}
\Lambda = \phi(x,y,z) + t\,\psi(x,y,z),
\end{equation}
where $\phi(x,y,z)$ and $\psi(x,y,z)$ must satisfy
\begin{equation}
\label{cond-gen_multipliers_KZK} 
\phi_{xx} + \phi_{yy} = \tfrac{2}{c_0} \psi_z, \quad
\psi_{xx} + \psi_{yy} = 0.
\end{equation}
The Laplace equation has infinitely many solutions $\psi(x,y;z)$,
called harmonic functions, which in this case are parameterized by
$z$. Once a solution for $\psi(x,y;z)$ is selected, one has to solve
the above Poisson equation for $\phi(x,y,z)$, using, e.g., the
Green's function approach. Since there are infinitely many solutions
for $\phi$ (with corresponding $\psi$) their exist infinite many
conservation laws for (\ref{KZK-org}). From
(\ref{characteristics_equation_KZK}), (\ref{gen_multipliers_KZK}),
and (\ref{cond-gen_multipliers_KZK}) one obtains
\begin{eqnarray}
\label{CL_gen_char_KZK}
&& \!\!\!\!\!\!\!\!\!\!\!\!\!\!\!\!\!\!
{\mathrm D}_{t}
\left[ \Lambda(\delta p_{tt} + 2 {\tilde{\beta}} p p_t-2c_0^3p_z)-
\Lambda_t(\delta p_t + {\tilde{\beta}} p^2) \right]
 + {\mathrm D}_{x}\, \left[ c_0^4 (\Lambda p_x-p\Lambda_x) \right]
\nonumber \\
&& \!\!\!\!\!\!\!\!\!\!\!\!\!\!\!\!\!\!
 + {\mathrm D}_{y} \left[ c_0^4 (\Lambda p_y-p\Lambda_y) \right]
 + {\mathrm D}_{z} \left[ 2c_0^3p \Lambda_t \right]
 \!=\! \Lambda \left(\delta p_{ttt} + {\tilde{\beta}} (p^{2})_{tt}
 - 2 c_0^3 p_{zt} + c_{0}^4 (p_{xx}+p_{yy}) \right)
  \\
&& \!\!\!\!\!\!\!\!\!\!\!\!\!\!\!\!\!\!
  - c_0^4 p (\Lambda_{xx} + \Lambda_{yy} - \tfrac{2}{c_0} \Lambda_{tz})
\nonumber
\end{eqnarray}
for arbitrary $p(x,y,z,t)$. Note that the last term in
(\ref{CL_gen_char_KZK}) vanishing since the last equation in
(\ref{cond-gen_multipliers_KZK}) holds. Then, when $p(x,y,z,t)$ is a
solution system (\ref{KZK}),
\begin{eqnarray}
\label{CL_KZK} &&
{\mathrm D}_{t} \, \left[\Lambda(\delta p_{tt}
  + 2 {\tilde{\beta}} p p_t-2c_0^3p_z)
  - \Lambda_t(\delta p_t + {\tilde{\beta}} p^2)\right]
  + {\mathrm D}_{x} \, \left[ c_0^4(\Lambda p_x-p\Lambda_x) \right]
\nonumber \\
&&
+ \, {\mathrm D}_{y} \, \left[ c_0^4(\Lambda p_y-p\Lambda_y)
\right] + {\mathrm D}_{z} \, \left[ 2c_0^3p\Lambda_t \right] \,
\dot{=} \, 0,
\end{eqnarray}
a local conservation law follows.
The conserved densities and fluxes for system (\ref{KZK})
with multiplier (\ref{gen_multipliers_KZK}) are
\begin{eqnarray}
\label{CL_KZK-densities-fluxes}
T^t &=& (\phi+t\psi)(\delta p_{tt} +2 {\tilde{\beta}} p p_t-2c_0^3p_z)
      - \psi (\delta p_t + {\tilde{\beta}} p^2),
\nonumber\\
T^x &=& c_0^4\left[(\phi+t\psi) p_x - p(\phi_x+t\psi_x)\right],
\\
T^y &=& c_0^4\left[(\phi+t\psi) p_y - p(\phi_y+t\psi_y)\right],
\nonumber \\
T^z &=& 2c_0^3p \psi,
\nonumber
\end{eqnarray}
where $\phi(x,y,z)$ and $\psi(x,y,z)$ are solutions of
(\ref{cond-gen_multipliers_KZK}).

\subsection{Conservation laws of the Khokhlov-Zabolotskaya-Kuznetsov equation
in cylindrical and spherical coordinates}
\label{KZK-equation-cylindrical}

Using the multiplier method, we
compute some conservation laws of the KZK equation
\cite[Ch.\ 8, Eq.\ (8.1)]{Hamilton-ch8-springer-2024}
in cylindrical and spherical coordinates,
\begin{equation}
\label{KZK_cylinder-original}
p_{rr} + \tfrac{p_r}{r} + \tfrac{\delta}{c_{0}^4} p_{ttt}
 + \tfrac{\beta}{\rho_0 c_{0}^4} (p^{2})_{tt} -\tfrac{2}{c_0} p_{zt} = 0,
\end{equation}
where the meaning of the symbols is analogous to (\ref{KZK-org}).
Taking a more general form,
\begin{equation}
\label{KZK_cylinder}
c_{0}^4 (p_{rr} + \tfrac{2 m p_r}{r}) + \delta p_{ttt}
+ {\tilde{\beta}} F(p)_{tt} - 2 c_0^3 p_{zt} = 0,
\end{equation}
with ${\tilde{\beta}} = \tfrac{\beta}{\rho_0}$, allows us to treat
the cylindrical ($m = \tfrac{1}{2}$) and spherical ($m = 1$) cases
and an arbitrary nonlinear function $F(p)$ at all once. Using
\verb|GeM| and the procedure outlined in the previous section then
yields the determining equations for the multiplier
$\Lambda(r,z,t,p, p_r,p_z,p_t)$:
\begin{equation}
\label{overdetermined_system_multipliers_KZK_general}
\Lambda_{tt} = 0, \;
\Lambda_{p} = 0, \;
\Lambda_{p_r} = 0, \;
\Lambda_{p_z} = 0, \;
\Lambda_{p_t} = 0, \;
\Lambda_{rr} = \tfrac{2}{c_0} \Lambda_{tz} + \tfrac{2m}{r^2}(r\Lambda_r-\Lambda),
\end{equation}
provided $F^{\prime\prime} \not = 0$. Solving
(\ref{overdetermined_system_multipliers_KZK_general}) one obtains
\begin{equation}
\label{multiplier_KZK_cylinder}
\Lambda = \phi(r,z) + t \, \psi(r,z),
\end{equation}
where $\phi(r,z)$ and $\psi(r,z)$ must be solutions of
\begin{equation}
\label{KZK-conditions-phi-psi}
\phi_{rr} -\tfrac{2 m}{r} \phi_{r} +\tfrac{2 m}{r^2} \phi
= \tfrac{2}{c_0} \psi_z,
\quad
\psi_{rr} -\tfrac{2 m}{r} \psi_{r} +\tfrac{2 m}{r^2} \psi = 0.
\end{equation}
The conserved density and flux for (\ref{KZK_cylinder}) with
multiplier $\Lambda(r,z,t,p,p_r,p_z,p_t)$ in
(\ref{multiplier_KZK_cylinder}) are
\begin{eqnarray}
\label{CL_KZK-density-flux}
T^t &=&
(\phi + t \psi)( \delta p_{tt} + {\tilde{\beta}} F'(p) p_t - 2 c_0^3 p_z )
      - \psi (\delta p_t + {\tilde{\beta}} F(p)),
\nonumber \\
T^r &=& c_0^4 \left( (\phi + t\psi) (p_r + \tfrac{2 m p}{r})
      - p (\phi_r+t\psi_r) \right),
\nonumber \\
T^z &=& 2 c_0^3 p \psi,
\end{eqnarray}
where $\phi(r,z)$ and $\psi(r,z)$ are solutions of
(\ref{KZK-conditions-phi-psi}).

\section{Application 2: Westervelt-type equations}
\label{application-Westervelt} Like the KZK equation, the Westervelt
equation is one of the most important nonlinear PDEs describing the
propagation of nonlinear waves in acoustics. We first compute
conservation laws for the Westervelt equation in one and two
dimensions before generalizing the results to any number of spatial
variables. In addition, we derive conservation laws for the
Westervelt equation in cylindrical and spherical coordinates.

\subsection{Conservation laws of Westervelt-type equations in Cartesian coordinates}
\label{Westervelt-equation-2D-3D}
\subsubsection{Conservation laws of Westervelt-type equations in two dimensions}
\label{Westervelt-equation-2D}

The Westervelt equation \cite[Ch.\ 3,
Eq.\ (3.46)]{Hamilton-Morfey-ch3-springer-2024}
\begin{equation}
\label{westervelt-2D-general}
\square^2 p + \tfrac{\delta}{c_0^4} p_{ttt}
+ \tfrac{\gamma}{\rho_0 c_0^4} (p^2)_{tt} = 0
\end{equation}
describes the propagation of acoustic waves with acoustic pressure
$p$ in various media. In Cartesian coordinates, $\square^2 =
\nabla^2 - \tfrac{1}{c_0^2} \frac{\partial^2}{\partial t^2}$ and
$\nabla^2 = \Delta = \frac{\partial^2}{\partial x^2} +
\frac{\partial^2}{\partial y^2}$ are the d'Alembertian and Laplacian
operators, respectively. In (\ref{westervelt-2D-general}), $\rho_0$
denotes the equilibrium density of the medium and $c_0$ is the speed
of sound; $\delta >0 $ is a diffusivity parameter (e.g., dissipation
or loss of sound pressure due to viscosity and thermal conduction in
the medium), and $\gamma >0 $ is a coefficient of nonlinearity.

Instead of analyzing (\ref{westervelt-2D-general}) we will work with
\begin{equation}
\label{Westervelt_2D}
F(p)_{tt} - c^2 (p_{xx} + p_{yy}) - \alpha p_{ttt}
  - \beta (p_{xx}+p_{yy})_t = 0,
\end{equation}
which includes the term in $\Delta_t \, p = \Delta p_t$ considered
in \cite[Eq.\ (10)]{Kaltenbacher-EECT-2025} with a
diffusivity-stabilizing parameter $\beta$, and an arbitrary
nonlinear function $F(p)$. For completeness we kept the third-order
term in time from (\ref{westervelt-2D-general}) which also appears
in some extended models discussed in \cite[Eqs.\ (63) and
(76)]{Kaltenbacher-EECT-2025} and references therein.

For $F(p) = \rho_0 c_0^2 p - \gamma p^2$, $c^2 = \rho_0 c_0^4$,
$\alpha = \delta \rho_0$, and $\beta = 0$, (\ref{Westervelt_2D})
becomes (\ref{westervelt-2D-general}) in two spacial dimensions. The
one-dimensional version of (\ref{westervelt-2D-general}) was
investigated by Anco et al.\
\cite{Anco-Marquez-Garrido-Gandarias-2023}. The conservation laws
obtained by these authors are a special case of those discussed in
Section~\ref{Westervelt-equation-nD}.

Equation (\ref{Westervelt_2D}) in one spatial dimension with $\alpha
= 0$ was studied by M\'{a}rquez et al.\ \cite[Eq.\
(3)]{Marquez-Recio-Gandarias-2025}. It arises in nonlinear
elasticity where $p$ then represents stress $\sigma$. See,
\cite[Eq.\ (1.1)]{Erbay-Sengul-2015} and \cite[Eqs.\ (22) and
(24)]{Sengul-DCDS-2021} with $F(p)$ as in Table 2.

Using \verb|GeM|, the determining equations for the multiplier
$\Lambda(x,y,t,p,p_x,p_y,p_t)$ are
\begin{equation}
\label{system_multipliers_Westervelt_2D}
\Lambda_{tt} = 0, \;\,
\Lambda_{p} = 0, \;\,
\Lambda_{p_x} = 0,\;\,
\Lambda_{p_y} = 0,\;\,
\Lambda_{p_t} = 0,\;\,
\Lambda_{xx} + \Lambda_{yy} = 0
\end{equation}
provided $F^{\prime\prime} \not = 0$. Solving
(\ref{system_multipliers_Westervelt_2D}) one gets
\begin{equation}
\label{multiplier_Westervelt_2D}
\Lambda = \phi(x,y) + t \, \psi(x,y),
\end{equation}
where $\phi(x,y)$ and $\psi(x,y)$ must both satisfy the Laplace
equation:
\begin{equation}
\label{conditions_Westervelt_2D} \phi_{xx} + \phi_{yy}  = 0, \quad
\psi_{xx} + \psi_{yy} = 0.
\end{equation}
The conservation law in characteristic form for the multipliers
(\ref{multiplier_Westervelt_2D}) satisfying
(\ref{conditions_Westervelt_2D}) reads
\begin{eqnarray}
\label{characteristics_equation_Westervelt_2D}
&& {\mathrm D}_{t} \,
T^t + {\mathrm D}_{x} \, T^x + {\mathrm D}_{y} \, T^y = (\phi(x,y) +
t \psi(x,y)) \left( F(p)_{tt} - c^2 (p_{xx} + p_{yy}) \right.
\nonumber \\
&& \left. - \alpha p_{ttt} - \beta (p_{xx}+p_{yy})_t \right),
\end{eqnarray}
for arbitrary $p(x,y,t)$. By integration of
(\ref{characteristics_equation_Westervelt_2D}), the density and flux
can be readily derived, resulting in
\begin{eqnarray}
\label{CL_char_Westervelt}
&& \!\!\!\!\!\!\!
{\mathrm D}_{t}
\left[(\phi+t\psi) \left(F^{\prime} p_t-\alpha p_{tt}\right)
- \psi (F(p)-\alpha p_t) \right]
+ {\mathrm D}_{x}
  \left[ -c^2 \left[(\phi+t\psi) p_x - p (\phi_x+t\psi_x)\right] \right.
\nonumber \\
&& \!\!\!\!\!\!\!
\left. -\beta [ (\phi+t\psi)p_{tx} - (\phi_x+t\psi_x) p_t ] \right]
+ {\mathrm D}_{y} \,
  \left[ -c^2 \left[(\phi+t\psi) p_y
   - p (\phi_y+t\psi_y) \right]
   - \beta \left[ (\phi+t\psi) p_{ty} \right.
  \right.
\nonumber\\
&& \!\!\!\!\!\!\!
 \left. \left. - (\phi_y+t\psi_y) p_t \right] \right]
= (\phi(x,y)+t\psi(x,y))
 \left( F(p)_{tt} - c^2 (p_{xx}+p_{yy}) - \alpha p_{ttt}
 - \beta (p_{xx}+p_{yy})_t \right)
\nonumber \\
&& \!\!\!\!\!\!\!
       + ( p c^2 +\beta p_t)
       \left( \phi_{xx}+\phi_{yy} + t (\psi_{xx}+\psi_{yy}) \right)
\end{eqnarray}
for arbitrary $p(x,y,t)$.
Note that the last term vanishes as a consequence of (\ref{conditions_Westervelt_2D}).

Evaluated on solutions $p(x,y,t)$ of (\ref{Westervelt_2D}) one then gets
\begin{eqnarray}
\label{CL_Westervelt_2D}
&& {\mathrm D}_{t} \,
   \left[ (\phi+t\psi) \left(F^{\prime} p_t-\alpha p_{tt} \right)
  - \psi( F(p)-\alpha p_t) \right]
  + {\mathrm D}_{x} \,
    \left[ -c^2 \left[(\phi+t\psi)p_x - p (\phi_x+t\psi_x) \right] \right.
\nonumber \\
&& \left.
  - \beta [ (\phi+t\psi)p_{tx} - (\phi_x+t\psi_x) p_t ] \right]
  + {\mathrm D}_{y} \, \left[-c^2 \left[(\phi+t\psi) p_y
  - p (\phi_y+t\psi_y) \right] \right.
  \nonumber \\
&& \left.
   - \beta [ (\phi+t\psi) p_{ty} -(\phi_y+t\psi_y) p_t ] \right]
   \, \dot{=} \, 0.
\end{eqnarray}
Therefore, the density and flux for (\ref{Westervelt_2D}) corresponding
to the multiplier family $\Lambda(x,y,t,p)$ in
(\ref{multiplier_Westervelt_2D}) are
\begin{eqnarray}
\label{CL_Westervelt_2D-density-flux}
T^t &=& (\phi+t\psi) \left( F^{\prime} p_t-\alpha p_{tt}\right)
        - \psi( F(p)-\alpha p_t),
\nonumber \\
T^x &=& -c^2 \left[(\phi+t\psi) p_x - p (\phi_x+t\psi_x) \right]
    - \beta [ (\phi+t\psi) p_{tx} - (\phi_x+t\psi_x) p_t ],
\\
T^y &=& -c^2 \left[(\phi+t\psi) p_y  - p (\phi_y+t\psi_y) \right]
        - \beta [ (\phi+t\psi) p_{ty} - (\phi_y+t\psi_y) p_t ],
\nonumber
\end{eqnarray}
where $\phi(x,y)$ and $\psi(x,y)$ are solutions of
(\ref{conditions_Westervelt_2D}).

\subsubsection{Conservation laws of Westervelt-type equations in three dimensions}
\label{Westervelt-equation-3D}

The three-dimensional version of the
Westervelt equation (in one spatial dimension) investigated by
M\'{a}rquez et al.\ \cite{Marquez-Recio-Gandarias-2025} reads
\begin{equation}
\label{Westervelt_3D}
F(p)_{tt} - \alpha p_{ttt} - \beta (p_{xx}+p_{yy}+p_{zz})_t
= c^2( p_{xx} + p_{yy} + p_{zz} ).
\end{equation}
Again using \verb|GeM|, the determining equations for
$\Lambda(x,y,z,t,p,p_x,p_y,p_z,p_t)$ are
\begin{equation}
\label{system_multipliers_Westervelt_3D}
 \Lambda_{tt} = 0, \;\,
 \Lambda_{p} = 0, \;\,
 \Lambda_{p_x} = 0, \;\,
 \Lambda_{p_y} = 0, \;\,
 \Lambda_{p_z} = 0, \;\,
 \Lambda_{p_t} = 0, \;\,
 \Lambda_{xx} + \Lambda_{yy} + \Lambda_{zz} = 0,
\end{equation}
provided $F^{\prime\prime} \not = 0$. Solving
(\ref{system_multipliers_Westervelt_3D}) yields
\begin{equation}
\label{multiplier_Westervelt_3D}
\Lambda = \phi(x,y,z) + t \, \psi(x,y,z),
\end{equation}
where $\phi(x,y,z)$ and $\psi(x,y,z)$ must solve the Laplace equation,
\begin{equation}
\label{conditions_Westervelt_3D} \phi_{xx} + \phi_{yy} + \phi_{zz} =
0, \quad
\psi_{xx} + \psi_{yy} + \psi_{zz} = 0.
\end{equation}
With the family of multipliers $\Lambda$ in
(\ref{multiplier_Westervelt_3D}), the density and flux for
(\ref{Westervelt_3D}) read
\begin{eqnarray}
\label{CL_Westervelt_3D-density-flux}
T^t &=& (\phi+t\psi)\left(F^{\prime} p_t-\alpha p_{tt}\right)
        - \psi(F(p)-\alpha p_t),
\nonumber \\
T^x &=& - c^2\left[(\phi+t\psi)p_x - p (\phi_x+t\psi_x) \right]
        - \beta [ (\phi+t\psi)p_{tx} - (\phi_x+t\psi_x)p_t ],
\\
T^y &=& -c^2\left[(\phi+t\psi)p_y - p (\phi_y+t\psi_y) \right]
        - \beta [ (\phi+t\psi)p_{ty} - (\phi_y+t\psi_y)p_t ],
\nonumber
\\
T^z &=& -c^2 \left[(\phi+t\psi)p_z - p (\phi_z+t\psi_z) \right]
         - \beta [ ( \phi+t\psi) p_{tz} - (\phi_z+t\psi_z) p_t ],
\nonumber
\end{eqnarray}
where $\phi(x,y,z)$ and $\psi(x,y,z)$ are any solutions of
(\ref{conditions_Westervelt_3D}).

\subsubsection{Conservation laws of Westervelt-type equations in any number of spatial  variables}
\label{Westervelt-equation-nD}

The Westervelt equation in $n$
independent variables $(x_1, x_2, \cdots, x_i, \cdots, x_n)$ denoted
below as $(x_i)$, and time $(t)$ is
\begin{equation}
\label{Westervelt_nD}
F(p)_{tt} - \alpha p_{ttt} - \beta \nabla^2 p_t = c^2 \nabla^2 p,
\end{equation}
where
\begin{equation}
\label{laplacian-nD}
\nabla^2
= \frac{\partial^2}{\partial x_1^2}
  + \frac{\partial^2}{\partial x_2^2}
  + \cdots
  + \frac{\partial^2}{\partial x_i^2}
  + \cdots
  + \frac{\partial^2}{\partial x_n^2}.
\end{equation}
Analogous to the two- and three-dimensional cases, the density and
flux for (\ref{Westervelt_nD}) with the family of multipliers of
type
\begin{equation}
\label{multiplier_Westervelt_nD}
\Lambda(x_i,t,p,p_{x_i},p_t) = \phi(x^i) +  t \, \psi(x^i)
\end{equation}
are
\begin{eqnarray}
\label{CLs_Westervelt_nD}
T^t &=& (\phi+t\psi) \left(F^{\prime} p_t-\alpha p_{tt}\right)
        - \psi (F(p) -\alpha p_t),
\nonumber \\
T^{x^i} &=&
 -c^2 \left[ (\phi + t\psi) p_{x^i} - p (\phi_{x^i} + t\psi_{x^i}) \right]
  - \beta [ (\phi + t\psi) p_{tx^i} - (\phi_{x^i} + t\psi_{x^i}) p_t ],
\end{eqnarray}
where $\phi(x^i)$ and $\psi(x^i)$ must be harmonic functions (i.e.,
solutions of the $n-$dimensional Laplace equation),
\begin{equation}
\label{conditions_Westervelt_ndim}
\nabla^2 \phi = 0, \quad
\nabla^2 \psi = 0.
\end{equation}
Using a theorem by Igonin \cite{Igonin-JPA-2024} and the direct
method \cite{Anco-Bluman-EurApplMath-I-2002,Olver-1993}, Sergyeyev
\cite{Sergyeyev-QTDS-2024} has derived the local conservation laws
of all orders of (\ref{Westervelt_nD}) with $\alpha \neq 0$ and
$\beta = 0$ (i.e., the lossless case). He has also shown that for
the one-dimensional Westervelt equation there are no other local
conservation laws than those reported by Anco et al.\
\cite[Eqs.\ (20)-(23)]{Anco-Marquez-Garrido-Gandarias-2023}.
In the one-dimensional case, apart from their linearity in $t$ expressed in
(\ref{multiplier_Westervelt_nD}), the functions $\phi(x)$ and
$\psi(x)$ are linear in $x$ because the conditions
(\ref{conditions_Westervelt_ndim}) reduce to $\phi_{xx} = \psi_{xx}
= 0$. Sergyeyev also gives a general formula for the conservation
laws of all orders (see, \cite[Eq.\ (2)]{Sergyeyev-QTDS-2024}) which
matches (\ref{CLs_Westervelt_nD}) again $\beta = 0$ and $c=1$.
Remarkably, the one-dimensional Westervelt equation has a finite
number of conservation laws whereas in multiple space variables the
equation has infinitely many conservation laws.

\subsection{Conservation laws of Westervelt-type equations
in spherical and cylindrical coordinates}
\label{Westervelt-equation-spherical}

In this section we consider
the Westervelt equation \cite[Eq.\ (7.1)]{Anco-EurApplMath-II-2023}
in a generalized form to cover the planar, cylindrical and spherical
cases,
\begin{equation}
\label{Westervelt_spherical-general}
p_{rr} + \tfrac{2 m}{r} p_r - p_{tt}  + \gamma (p^2)_{tt}
+ \alpha p_{ttt} + \beta(p_{rr} + \tfrac{2 m}{r} p_{r})_t = 0,
\end{equation}
where, for convenience, the units are chosen
\cite{Anco-Marquez-Garrido-Gandarias-2023} such that $c=1$.
Note that $m = 0$ corresponds to the planar case (replacing $r$ by $x$);
$m = \tfrac{1}{2}$ covers the cylindrical case; and $m=1$ the spherical case.
The latter case is treated in \cite[Eq.\ (7.1)]{Anco-EurApplMath-II-2023}.

For the planar case ($m = 0$ with $x$ instead of $r$),
(\ref{Westervelt_spherical-general}) can be written as
\begin{equation}
\label{Westervelt_planar}
\big( (p - \gamma p^2 - \alpha p_t)_t - \beta p_{xx} \big)_t + (-p_x)_x = 0,
\end{equation}
which is a conservation law itself.

Using \verb|GeM| and a multiplier of type $\Lambda(r,t,p,p_r,p_t)$,
we computed the conservation laws for
(\ref{Westervelt_spherical-general}) presented in Table 5 for
arbitrary values of $m$ as well as for $m=1$ and $m=\tfrac{1}{2}$.
\begin{center}
\begin{table}[htb]
\label{table-CLs-Westervelt-spherical}
\small
\begin{tabular}{llll}
\hline $m$ & $\Lambda$ & Density $T^r$ & Flux $T^t$
\\
\hline
\\[-0.5em]
\\[-0.5em]
& $ r^{2m}$ & $ - r^{2m}( p_r+\beta p_{rt})$ &
$ r^{2m} \big((1-2\gamma p) p_t- \alpha p_{tt} \big)$
\\[-0.5em]
\\[-0.5em]
& $ r$ & $ (1-2 m) p - r p_r $ &
$r \big( (1-2 \gamma p) p_t -\alpha p_{tt} \big)$
\\[-0.7em]
\\
& & $ + \beta \big( (1-2m) p_t-rp_{tr} \big)$ &
\\[-0.5em]
\\
& $ r^{2m} t$ & $-t r^{2m} (p_r+\beta p_{rt})$ & $r^{2m}
\big( t \big( (1-2\gamma p)p_t - \alpha p_{tt} \big)
  - (1-\gamma p) p + \alpha p_t \big)$
\\[-0.5em]
\\
& $ r t$ & $ t \big( (1-2m)p -r p_r \big)$ & $r \big( t \big(
(1-2\gamma p)p_t - \alpha p_{tt} \big)
  - (1-\gamma p) p + \alpha p_t \big)$
\\[-0.5em]
\\
& & $ +\beta t \big( (1-2m)p_t-rp_{rt} \big)$ &
\\[-0.5em]
\\
\hline
\\[-0.5em]
\\[-0.5em]
$1$ & $ r^{2}$ &  $ - r^{2}( p_r+\beta p_{rt})$ & $ r^{2} \big(
(1-2\gamma p) p_t- \alpha p_{tt} \big)$
\\[-0.5em]
\\
& $ r$ & $ -(p + r p_r)-\beta(p_t+rp_{rt})$ & $
r \big( (1-2 \gamma p) p_t -\alpha p_{tt} \big)$
\\[-0.5em]
\\
& $ r^{2}t$ & $- r^{2} t (p_r+\beta p_{tr})$ &
$ r^{2} \big( t \big( (1-2\gamma p)p_t - \alpha p_{tt} \big)
  - (1-\gamma p) p + \alpha p_t \big)$
\\[-0.5em]
\\
& $ r t$ & $ -t \big( p +r p_r \big) - \beta t(p_t+rp_{rt})$ &
$ r \big( t \big( (1-2\gamma p)p_t - \alpha p_{tt} \big)
  - (1-\gamma p) p + \alpha p_t \big)$
\\[-0.5em]
\\
\hline
\\[-0.5em]
\\[-0.5em]
$ \frac{1}{2}$ & $ r$ & $ -r ( p_r+\beta p_{rt})$ & $ r \big( (1-2
\gamma p) p_t -\alpha p_{tt} \big)$
\\[-0.5em]
\\
& $r t$ & $- r t( p_r+\beta p_{rt})$ & $
r \big( t \big( (1-2\gamma p) p_t - \alpha p_{tt} \big)
  - (1-\gamma p) p + \alpha p_t \big)$
\\[-0.5em]
\\
& $r \ln r$ & $p - r p_r \ln r +\beta(p_t-r \ln(r) p_{tr})$ &
$ r \ln(r) \big( (1-2 \gamma p) p_t -\alpha p_{tt} \big)$
\\[-0.5em]
\\
& $r t \ln r$ & $t \big(p - r p_r \ln r\big) $ &
$ r \ln(r) \big( t \big( (1-2\gamma p)p_t - \alpha p_{tt} \big)
   - (1-\gamma p) p + \alpha p_t \big)$
\\[-0.6em]
\\
& & $ +\beta t \big( p_t-r\ln(r) p_{tr} \big)$ &
\\[-0.5em]
\\
\hline
\end{tabular}
\caption{
Multipliers, conserved densities, and fluxes for
(\ref{Westervelt_spherical-general}) for arbitrary $m$ (top),
$m=1$ (middle), and $m = \tfrac{1}{2}$ (bottom).
}
\end{table}
\end{center}

The conservation laws for the one-dimensional Westervelt equation
presented in \cite[Eq.\ (1)]{Anco-Marquez-Garrido-Gandarias-2023}
can be obtained by setting $m = \beta = 0$ and replacing $T^r$ by
$T^x$ in Table 5, yielding the four polynomial conservation laws
derived by Anco et al.\ \cite[Eqs.\
(20)-(23)]{Anco-Marquez-Garrido-Gandarias-2023}.
For $m = \alpha = \beta = 0$ there also exists a rational conservation law
(not included in Table 5) given in
\cite[Eq.\ (24)]{Anco-Marquez-Garrido-Gandarias-2023}.
For the spherical case $(m=1)$ and $\beta = 0$ the first two fluxes
$(T^t)$ in Table 5 correspond to the conserved densities given
in \cite[Eqs.\ (7.8) and (7.9)]{Anco-EurApplMath-II-2023}.
For the cylindrical case ($m=\tfrac{1}{2}$) the first and second
conservation laws for arbitrary $m$ in Table 5 coincide and so do
the third and fourth ones.

Next, replacing the nonlinearity $p^2$ in (\ref{Westervelt_spherical-general})
by an arbitrary function $F(p)$, we compute conservation laws for
\begin{equation}
\label{Westervelt_spherical-general_F}
p_{rr} + \tfrac{2 m}{r} p_r + F(p)_{tt} + \alpha p_{ttt}
+ \beta ( p_{rr} + \tfrac{2 m}{r} p_{r} )_t = 0,
\end{equation}
where, as above, $m = 0, \tfrac{1}{2}$, or $1$. The results are
presented in Table 6.

\begin{center}
\begin{table}[htb]
\label{table-CLs-Westervelt-spherical_F}
\small
\begin{tabular}{llll}
\hline $m$ & $\Lambda$ & Density $T^r$ & Flux $T^t$
\\
\hline
\\[-0.5em]
\\
& $ r^{2m}$ & $ - r^{2m}( p_r+\beta p_{rt})$ &
$ -r^{2m} ( F^{\prime} p_t + \alpha p_{tt} )$
\\[-0.5em]
\\
& $ r$ & $ (1-2 m) p - r p_r + \beta \big((1-2m)p_t-rp_{tr} \big)$ &
$ -r ( F^{\prime} p_t + \alpha p_{tt} )$
\\[-0.5em]
\\
& $ r^{2m}t$ & $-t r^{2m} (p_r+\beta p_{rt})$ &
$ -r^{2m} \big( t ( F^{\prime} p_t + \alpha p_{tt} )
   - F(p) - \alpha p_t \big)$
\\[-0.5em]
\\
& $ r t$ & $
t \big( (1-2m)p - r p_r \big) + \beta t \big( (1-2m)p_t-rp_{rt} \big) $ &
$ - r \big( t ( F^{\prime} p_t + \alpha p_{tt} )
  - F(p) - \alpha p_t \big)$
\\[-0.5em]
\\
\hline
\\[-0.5em]
\\
$1$ & $ r^{2}$ &  $ - r^{2}( p_r+\beta p_{rt})$ & $ - r^{2} (
F^{\prime} p_t + \alpha p_{tt} )$
\\[-0.5em]
\\
& $ r$ & $ -(p + r p_r)-\beta(p_t + r p_{rt})$ &
$ -r ( F^{\prime} p_t + \alpha p_{tt} )$
\\[-0.5em]
\\
& $ r^{2} t$ & $- r^{2} t (p_r+\beta p_{tr})$ & $
-r^{2} \big( t ( F^{\prime} p_t + \alpha p_{tt} )
   - F(p) - \alpha p_t \big)$
\\[-0.5em]
\\
& $ r t$ & $ -t \big( p +r p_r \big)-\beta t(p_t+rp_{rt})$ &
$ - r \big( t ( F^{\prime} p_t + \alpha p_{tt} )
  - F(p) - \alpha p_t \big)$
\\[-0.5em]
\\
\hline
\\[-0.5em]
\\
$ \frac{1}{2}$ & $ r$ & $ -r ( p_r+\beta p_{rt})$ & $ - r (
F^{\prime} p_t + \alpha p_{tt} )$
\\[-0.5em]
\\
& $r t$ & $- r t( p_r+\beta p_{rt})$ &
$ - r \big( t ( F^{\prime} p_t + \alpha p_{tt} )
   - F(p) - \alpha p_t \big)$
\\[-0.5em]
\\
& $r \ln r$ & $p - r p_r \ln r + \beta(p_t-r \ln(r) p_{tr})$ &
$ - r \ln(r) ( F^{\prime} p_t + \alpha p_{tt} ) $
\\[-0.5em]
\\
& $r t \ln r$ &
$t \big( p - r p_r \ln r +\beta(p_t-r \ln(r) p_{tr}) \big) $ &
$ - r \ln(r) \big( t ( F^{\prime} p_t + \alpha p_{tt} )
   - F(p) - \alpha p_t \big)$
\\[-0.5em]
\\
\hline
\end{tabular}
\caption{
Multipliers and conserved densities and fluxes for
(\ref{Westervelt_spherical-general_F}) for arbitrary $m$,
and $m = 1$ and $m = \tfrac{1}{2}$.
}
\end{table}
\end{center}

The conservation laws for the one-dimensional Westervelt equation in
\cite[Eq.\ (3)]{Marquez-Recio-Gandarias-2025} are obtained by
setting $m = \alpha = 0$, $c = 1$, $F(p) = -f(p)$, and replacing $r$
by $x$ in Table 6. Indeed, the first four conservation laws then
correspond to those computed by M\'{a}rquez et al.\ \cite[Eqs.\
(21)-(24)]{Marquez-Recio-Gandarias-2025}.

\section{Conclusions}
\label{conclusions}
In this paper, conservation laws have been
computed for several nonlinear PDEs arising in elasticity and
acoustics using two methods: the scaling-homogeneity approach which
is implemented in the {\emph{Mathematica}}-based code
\verb|ConservationLawsMD.m| and the multiplier method available in
the code \verb|GeM| written in {\emph{Maple}} syntax.

First, conservation laws have been derived for two models describing
shear wave propagation in a circular cylinder and a cylindrical annulus.
The PDEs are based on a constitutive relation expressed as a power law
for the stress with an arbitrary exponent $n$.
Computation of conservation laws with the scaling-homogeneity approach
required fixed values of $n$.
Based on the results for $n = 1, 2, \cdots$, the densities and fluxes
for arbitrary $n$ were readily obtained by pattern matching or
with a bit of extra work.
The computations are shown step-by-step and in sufficient detail
for the reader to be able to reproduce the results.

Second, with the aid of \verb|GeM|, the multiplier approach was applied
to a parameterized system involving an arbitrary function of stress.
The general system not only includes the two previous models
(for specific values of the parameters) but covers a broad class of models
for wave propagation in elastic materials formulated in terms of
cylindrical coordinates.

Next, conservation laws have been computed for the KZK equation in Cartesian,
cylindrically, and spherical coordinates.
For multipliers involving up to first-order derivatives of the dependent variables, \verb|GeM| allowed us to compute and simplify the determining equations for the multipliers.
These equations were then solved by hand showing that the multipliers are
linear in time $(t)$.
In Cartesian coordinates the space-dependent coefficient of $t$
satisfies the Laplace equation while the other function must be a solution
of Poisson's equation.
In cylindrical and spherical coordinates the situation is similar:
the multipliers are still linear in time but the equations for its coefficients
are more complicated.
Since the Laplace equation has infinitely many solutions
(called harmonic functions)
an infinite number of conservation laws for the KZK equation is guaranteed.
The closed form expressions for the conserved densities and fluxes have been
obtained by straightforward algebraic manipulations.

As a final application, conservation laws have been derived with the
multiplier method and \verb|GeM| for Westervelt-type equations in
multi-spatial (Cartesian) coordinates as well as cylindrical and
spherical coordinates. In all cases, the commonly used quadratic
nonlinearity was replaced by an arbitrary function of the acoustic
pressure. As with the KZK equation, the multipliers are linear in
time but in this case both coefficients must be solutions of the
Laplace equation. Therefore, the Westervelt-type equations in
Cartesian coordinates with more than one spatial variable have
infinitely many conservation laws which confirms recent findings by
Sergyeyev. By contrast, in one spatial variable the Westervelt-type
equation has only a finite number of conserved densities and fluxes
which have been reported in the literature.

\section*{Dedication}
\label{dedication} This paper is dedicated to Dr.\ George Bluman who
has been an inspiring teacher and leader in the symmetry community.
We are grateful for his seminal contributions to the advancement of
symmetry methods for differential equations. His journal articles
and books have greatly impacted the development of our own research.
We owe him our profound gratitude and wish him a long and happy
retirement.

\section*{Acknowledgments}
\label{acknowledgments}
We would like to thank Dr.\ Giuseppe Saccomandi
(University of Perugia, Perugia, Italy) for valuable
comments. We are very grateful to Dr.\ Roger Bustamante (University
of Chile, Santiago, Chile) for insightful comments and providing his
review papers on implicit constitutive theories and wave propagation
in elasticity bodies. We thank Dr.\ Norbert Euler for inviting us to
contribute to this special issue of OCNMP in honor of Dr.\ Bluman.

\label{lastpage}
\end{document}